\documentclass[12pt]{amsart}
\usepackage[letterpaper,margin=1in]{geometry}
\usepackage{graphicx,amsfonts,amsmath}
\usepackage{amssymb,amsthm,fancyhdr,indentfirst} 
\usepackage{amscd}
\usepackage{bm} 
\usepackage{color}
\usepackage{tikz}
\usetikzlibrary{shapes}
\usepackage[ruled,vlined]{algorithm2e}

\def\ff{\mathit{f\mkern-5mu f}}
\def\fc{\mathit{f\mkern-5mu c}}
\def\cf{\mathit{c\mkern-2mu f}}
\def\cc{\mathit{c\mkern-3mu c}}

\usepackage{color}
\theoremstyle{definition} 
\theoremstyle{plain} 
\theoremstyle{plain} 
\theoremstyle{plain} 
\theoremstyle{plain} 
\theoremstyle{plain} \newtheorem{lemma}{Lemma}
\theoremstyle{plain} \newtheorem{corollary}{Corollary}
\theoremstyle{plain} \newtheorem{remark}{Remark}
\theoremstyle{plain} \swapnumbers 
\numberwithin{equation}{section}
\numberwithin{figure}{section}
\numberwithin{table}{section}
\newcommand{\tvV}{\mathcal{V}}

\def\XNorm#1{\left\| #1 \right\|}                     
\def\XIProd#1#2{\left\langle #1 ,~ #2 \right\rangle}  
                                  
\def\XVec#1{{\mathbf #1}}       

\newcommand{\innerprod}[3][2]{\langle #2,#3 \rangle_{#1}}

\def\ZS{S}

\def\ZP{R}

\def\Xv{\XVec{e}}

\def\Xv{\XVec{v}}

\def\XM{\mu}                   
          
\def\XMx{\XM^\star}

\def\Xv{\XVec{e}}

\def\Xv{\XVec{v}}

\def\XMt{{\widetilde{M}}}

\newcommand{\norm}[2][2]{\| #2 \|_{#1}}

\newcounter{Xvqnctr}

\def\XMt{{\widetilde{M}}}
\newcommand{\be}{\begin{equation}}
\newcommand{\ee}{\end{equation}}

\newcommand{\SetAlgorithmStyle}{
	\SetKwData{Left}{left}\SetKwData{This}{this}\SetKwData{Up}{up}
	\SetKwInOut{Input}{Input}\SetKwInOut{Output}{Output}
	\ResetInOut{input}
	\SetKwComment{tcp}{//}{}
	\SetKwFor{For}{for}{}{end}
	\SetKwFor{While}{while}{}{endw}
	\SetArgSty{}
	\DontPrintSemicolon
}

\def\Popt{P_{\sharp}}      
\def\PoptP{\bar{ P_{\sharp} } }   
\def\RoptR{\bar{ R_{\sharp} } }   
  
\def\Pid{P_{\star}}      
\def\Zid{Z_{\star}}  	

\def\Sopt{S_{\sharp}} 

\def\PiMt{\Pi_{\XMt}} 
\begin{document}
\author{James Brannick, Fei Cao, Karsten Kahl, Rob Falgout and Xiaozhe Hu}
\title[Optimal interpolation and Compatible Relaxation in Classical AMG]{Optimal interpolation and Compatible Relaxation in Classical Algebraic Multigrid}
\begin{abstract}
In this paper,  we consider a classical form of  {\em optimal} algebraic multigrid (AMG) interpolation that directly minimizes 
the two-grid convergence rate and compare it with the so-called ideal form that
minimizes a certain {\em weak approximation property} of the coarse space.  
We study compatible relaxation type estimates for the quality of the coarse grid
and derive a new sharp measure using optimal interpolation that provides a guaranteed lower
bound on the convergence rate of the resulting two-grid method for a given grid.  In addition, we design a 
generalized bootstrap algebraic multigrid setup algorithm that computes
a sparse approximation to the optimal interpolation matrix.  We demonstrate numerically that
the BAMG method with sparse interpolation matrix (and spanning multiple levels) outperforms the two-grid method with
the standard ideal interpolation (a dense matrix) for various scalar 
diffusion problems with highly varying diffusion coefficient.
\end{abstract}
\maketitle

\section{Introduction}
We analyze algebraic multigrid (AMG) coarsening algorithms for linear systems of algebraic equations
\begin{equation}\label{eq:mod}
A \bf{u} = \bf{f},
\end{equation} 
coming from cell-centered finite volume discretizations of the scalar elliptic diffusion problem
\begin{equation} \label{eqn:model}
\begin{cases}
- \nabla \cdot (a(x) \nabla u)  = f \quad \text{in} \ \Omega \\
 u = 0 \qquad\qquad\quad\qquad  \text{on} \ \partial \Omega
\end{cases}
\end{equation}
where $\Omega = [0,1] \times [0,1]$.  
The discrete solution and right-hand side satisfy 
$\bf{u},\bf{f} \in \mathbb{R}^{n}$ and $A \in \mathbb{R}^{n\times n}$
is a symmetric and positive definite matrix.  We consider the case where the diffusion 
coefficient $a(x)$ is highly oscillatory, which is a problem that motivated the design of 
the original classical AMG setup algorithm~\cite{1985BrandtA_McCormickS_RugeJ-aa,1985BrandtA_McCormickS_RugeJ-ab}.

Multigrid solvers for solving~\eqref{eqn:model} involve a smoother, $M$, with error propagator given by $I - M^{-1}A$ and a coarse-level correction with error propagator given by $I - \Pi_A = I - PA_c^{-1}P^TA$, where
$P \in \mathbb{R}^{n \times n_c}$ denotes the interpolation matrix and $A_c = P^TAP$ is the Galerkin variational coarse-level matrix. The error propagation matrix of the resulting two-level method 
reads 
\begin{equation}\label{eq:tl}
E_{TG} = (I-M^{-1}A)(I-\Pi_A).
\end{equation} 
In AMG, the smoother $M$ is typically fixed and then interpolation $P$ is constructed in
an automated setup algorithm that takes as input the system matrix $A$ and computes
$P$ and $A_c$.
The main task in the AMG setup algorithm is thus to construct a {\em stable} interpolation matrix $P$ such 
that a certain {\em approximation property} holds and both $P$ and $A_c$ are {\em sparse} matrices.
The latter sparsity requirement implies that the procedure can be applied recursively in order to construct an optimal multilevel 
solver.

Numerous setup algorithms have been developed for constructing matrix-dependent 
interpolation, going back to the original classical AMG algorithm~\cite{1985BrandtA_McCormickS_RugeJ-aa,1985BrandtA_McCormickS_RugeJ-ab}. 
Generally speaking, the setup algorithm for constructing $P$ can be separated into three tasks:
\begin{enumerate}
\item Choosing the set of coarse variables, $C$, with cardinality $n_c= |C|$.
\item Determining the nonzero sparsity structure of $P$. 
\item Computing the values of the nonzero entries in $P$.
\end{enumerate}
Oftentimes, steps (1) and (2) of the setup algorithm are combined into a single step, 
as in smoothed aggregation AMG (SA)~\cite{PVanek_JMandel_MBrezina_1995a} where the choice of aggregates also determines 
the sparsity structure of the columns of $P$.  The coefficients of the columns of $P$ are then 
chosen to approximate certain error components that the smoother cannot treat efficiently, assumed in most cases
to be error that is dominated by the eigenvectors of the system matrix with small eigenvalues.  
In other approaches, e.g., classical AMG, steps (1)-(3)
are implemented in different stages within the setup algorithm~\cite{1985BrandtA_McCormickS_RugeJ-aa}.  
Specifically, the notion of strength of coupling between {\em neighboring} unknowns
is used in a maximal independent set algorithm to choose the coarse variable set $C$.  Next, the strongly and weakly
coupled neighbors of each of the fine degrees of freedom are determined and, finally, the coefficients 
of the corresponding row of interpolation are computed in a way that ensures that certain components
of the error (e.g., the constant vector) are well approximated locally.  
We note that in both approaches, once the coarse degrees of freedom $C$ and the entries of the 
interpolation matrix $P$ are selected they are fixed and the setup algorithm proceeds to construct the next coarser level.  
In this way, the algorithm is applied recursively without any measure of the quality of the resulting
coarse space in approximating the error in the current solution for the fine-level equations.

Compatible relaxation (CR)~\cite{PBrandt_2000,JBrannick_RFalgout,Falgout} and  adaptive~\cite{Breetal2,Breetal3} and bootstrap 
AMG~\cite{BAMG2,BAMG2010,BAMG14,BKL14,BKS14} were introduced as techniques to modify and adjust
the coarse variable set and interpolation, respectively.  The basic idea in these approaches is to  develop
local measures to assess the
suitability of the computed coarse space for a given problem.  Although these approaches have been successfully developed
and extended to handle numerous applications, certain theoretical issues remain unresolved.  For 
example, although the convergence rates of the CR algorithms that are typically used in practice give a qualitative measure of the suitability of the coarse set, they do
not accurately predict the convergence rate of the resulting two-grid solver for ideal interpolation in general 
(see~\cite{BBKL14b,JBrannick_RFalgout}), which is the aim of the approach.  Moreover, the so-called ideal interpolation matrix used as the basis of CR does not
in general give the fastest possible convergence rate of the two-level method over all possible choices of $P$ and, hence, it may 
not provide a reliable measure of the quality of the coarse variable set for certain problems.

In this paper, we study these issues further with the aim of gaining a deeper understanding of 
AMG from theoretical and practical points of view.
In Section~\ref{sec:prelim}, we derive the optimal classical AMG interpolation matrix and then
contrast it with the so-called ideal form that
minimizes a certain {\em weak approximation property} of the coarse space.
Section~\ref{sec:sharpCR}  introduces measures of the quality of the coarse grid based on the notions of
compatible relaxation~\cite{PBrandt_2000,JBrannick_RFalgout} and ideal interpolation as well as 
this new optimal form of classical AMG interpolation.  We show that the reliability and robustness 
of CR depends critically on the choice of the coarse variable type and that when the simplified (and computable) $F$-relaxation
form of CR is used, the resulting estimates of the two-grid solver with ideal interpolation and full smoothing are not sharp in general.  Then, 
we derive an iteration for accurately approximating the convergence of the two-grid method with ideal $P$ and show that it 
can be efficiently computed in practice.  We note, however, that even this more accurate CR-type estimate is not reliable in general since
for certain scalar diffusion test problems that we consider the optimal $P$ results in significantly faster convergence.  
To address this limitation, we derive a sharp variant of CR based on the optimal classical AMG interpolation matrix.  
On the other hand, we derive an equivalence between the optimal form of classical AMG interpolation and the so-called ideal interpolation
matrix in the case that $F$-smoothing is used in the resulting two-grid solver.   In addition, we derive a generalization of the ideal interpolation operator and show that for proper choices of the coarse variables this generalized ideal interpolation is equivalent to the optimal form.  Section \ref{sec:num} contains numerical results that illustrate these findings for
scalar diffusion test problems.  In addition, this section contains the derivation of a generalized bootstrap AMG (BAMG) setup algorithm
that aims to approximate the optimal interpolation matrix with a sparse approximation.  The main new feature of the algorithm is that it 
computes approximations to eigenvectors with small eigenvalues of the generalized eigenvalue problem for $(A,\XMt)$, where $A$ denotes 
the system matrix and $\XMt$ the symmetrized smoother.  Numerically, we show that the BAMG method (spanning multiple levels) with sparse $P$ outperforms the two-grid method with the ideal $P$ (which is a dense matrix) for our test problems. 

\section{Two-level theory and optimal classical AMG interpolation}\label{sec:prelim}

In~\cite{FVZ}, the following identity for convergence rate of the two-level method is introduced
\begin{equation}\label{eq:twogridbound}\displaystyle{
\|E_{  TG}(P)\|_A^2 = 1 - \dfrac{1}{\sup_{\Xv}\kappa(P,\Xv)}}, \quad
\kappa(P,\Xv) = \dfrac{\|(I - \PiMt(P))\Xv\|^2_{\XMt}}{\|\Xv\|^2_A},
\end{equation} where $\PiMt(P)$ is the $(\cdot,\cdot)_{\XMt}$ orthogonal projection on $\operatorname{range}(P)$,
with $\XMt^{-1}=M^{-1}+M^{-T}-M^{-T}AM^{-1}$ denoting the symmetrized smoother, so that $\XMt=M(M+M^T-A)^{-1}M^T$.
Note that, assuming $M+M^T-A$ is symmetric and positive definite (SPD) is equivalent to assuming the convergence of the chosen smoother defined by $M$.

Using~\eqref{eq:twogridbound} it is straightforward to derive the optimal two-grid convergence rate $\|E_{  TG}(P)\|_A^2$ with respect to $P$, for a given smoother $M$.  We note that this result is found in~\cite{FVZ} (see Corollary 4.1) and  more recently in the review paper~\cite{XZ_AMG}. 
From this general form of optimal $P$ we then derive
the classical AMG form of optimal P via a post-scaling, which is the result of primary interest in this paper.  In particular, we focus on studying this optimal classical AMG interpolation and comparing it in a practical setting with the so-called ideal classical AMG interpolation matrix. 

\begin{lemma}\label{lem:optimalP}
     Let $P: \mathbb{R}^{n_{c}} \rightarrow \mathbb{R}^{n}$ be full rank 
and let $\lambda_{1}\leq \lambda_{2} \leq \ldots \leq \lambda_{n}$ and 
$\Xv_{1},\Xv_{2},\ldots,\Xv_{n}$
     denote the eigenvalues and orthonormal (w.r.t $(\cdot,\cdot)_\XMt$ for convenience of representation )
eigenvectors of the generalized eigenvalue problem
     \begin{equation}\label{eq:gep}
         A{\bf x} = \lambda \XMt {\bf x}.
     \end{equation} Then the minimal convergence rate of the two-grid 
method is given by
     \begin{equation}\label{eq:sharp}
         \|E_{  TG}(\Popt)\|_A^2 = 1 - \kappa_{\sharp}, \quad \kappa_{\sharp} :=  \frac{1}{\displaystyle{\inf_{\operatorname{dim}(\operatorname{range}(P))\, =\, n_{c}} \sup_v \kappa(P,\Xv)}} =  \lambda_{n_{c}+1},
     \end{equation} where the optimal interpolation operator $\Popt$ 
satisfies
     \begin{equation}\label{eq:sharpP}
         \operatorname{range}(\Popt) = \operatorname{range}( 
\begin{pmatrix} \Xv_{1} & \Xv_{2} & \cdots & \Xv_{n_{c}} \end{pmatrix}).
     \end{equation}
     For sake of definiteness we set $\Popt = \begin{pmatrix} \Xv_{1} & 
\Xv_{2} & \cdots & \Xv_{n_{c}} \end{pmatrix}$ throughout the paper.
     \begin{proof} Starting with the equality in~\eqref{eq:twogridbound} 
first observe that if we denote $\mathcal{P} = \operatorname{range}(P)$ 
and its $\XMt$-orthogonal complement by $\mathcal{S}$ we find that for 
any $P$
         \begin{equation*}
             \dfrac{1}{\sup_{\Xv}\kappa(P,\Xv)} = \inf_{\Xv} 
\dfrac{\|\Xv\|^2_A}{\|(I - \pi_{\XMt}(P))\Xv\|^2_{\XMt}} \leq
\inf_{\Xv\, \in\, \mathcal{S}} \dfrac{\|\Xv\|_{A}^2}{\|\Xv\|_{\XMt}^2}.
         \end{equation*} Thus, we obtain for any $P$ that
         \begin{equation*}
             \|E_{  TG}(P)\|_A^2 \geq 1 - \inf_{\Xv\, \in\,
\mathcal{S}} \dfrac{\|\Xv\|_{A}^2}{\|\Xv\|_{\XMt}^2}
         \end{equation*} Finally due to 
$\operatorname{dim}(\mathcal{S}) + \operatorname{dim}(\mathcal{P}) = n$ 
we get
         \begin{equation*}
             \inf_{\operatorname{dim}(\mathcal{P})\, =\, n_{c}} \|E_{  
TG}(P)\|_{A}^{2} \geq
             1 - \sup_{\operatorname{dim}(\mathcal{S})\, =\, n - n_{c}} 
\inf_{\Xv\, \in\, \mathcal{S}\setminus 
\{0\}}\dfrac{\|\Xv\|^2_{A}}{\|\Xv\|^2_{\XMt}}.
         \end{equation*} 
Based on Courant-Fischer Min-max representation, we obtain:
         \begin{equation*}
             \inf_{\operatorname{dim}(\mathcal{P})\, =\, n_{c}} \|E_{  
TG}(P)\|_{A}^{2} \geq 1 - \lambda_{n_{c}+1}.
         \end{equation*}
The equality in this bound is obtained by setting $\mathcal{S} =\mathcal{S}_{\sharp}= \operatorname{range}(\begin{pmatrix} 
\Xv_{n_{c}+1} & \Xv_{n_{c}+2} & \cdots & \Xv_{n} \end{pmatrix})$.   
         Now, since $\mathcal{P}$ is the $\XMt$-orthogonal complement of 
$\mathcal{S}$, it follows that $P=\Popt$ and 
 \begin{equation*}
 \|E_{  TG}(\Popt)\|_{A}^{2} =1 - \lambda_{n_{c}+1}.
\end{equation*}
Finally, the optimal convergence rate for the two-grid method is obtained by choosing any interpolation that has the same range as $\Popt$, namely,
by setting $$P=\Popt Z, \quad \text{where we assume}\quad Z^{-1} \quad \text{exists.}$$   In this way, one obtains the same optimal convergence rate since by direct computation 
\begin{align*}
E_{TG}(\Popt Z) &= (I-M^{-1}A)(I-\Popt Z((\Popt Z)^TA\Popt Z)^{-1}(\Popt Z)^TA ) \\ 
&  = (I-M^{-1}A)(I-\Popt ((\Popt )^TA\Popt )^{-1}(\Popt )^TA )  =: E_{TG}(\Popt).
\end{align*} 
\end{proof}
\end{lemma}

We note that the above identity for the projection on $\Popt$ holds for $\Pi_X(\Popt)$ as well, where X is assumed to be
any SPD matrix.  This is summarized in the following corollary.
\begin{corollary}\label{eq:invariant}
Any projection $\Pi_X(P) = P(P^TXP)^{-1}P^TX$ is invariant with respect to post-multiplication of interpolation $P$ by an invertible matrix $Z$, $P \leftarrow PZ$.  Here $X$ is assumed to be
any SPD matrix.  
\begin{proof}
The proof is identical to the derivation for $E_{TG} (\Popt Z)$, where the system matrix $A$ is replaced by $X$ and the smoother is omitted.
\end{proof}
\end{corollary}

In order to derive the standard classical AMG form of the optimal interpolation $\Popt$, we use the fact that
the spectral radius of $E_{TG}$ remains unchanged if we replace $\Popt$ by $\Popt Z$ for any
nonsingular (invertible) matrix $Z$.   

\begin{remark}\label{cor:clOptP}
The classical AMG form of interpolation, assuming a splitting
of the fine-level degrees of freedom into $C$ and $F$, is given by
\begin{equation}\label{eq:classP}
P = \left[\begin{matrix}
   W \\ I 
\end{matrix} \right]
\begin{array}{lr}
\} \: F \\ \} \:  C
\end{array} ,
\end{equation}
where $C \subset  \{1,...,n\}$ and $F =  \{1,...,n\} \setminus C$, and $W \in \mathbb{R}^{|F| \times |C|}$ 
defines the interpolation weights.  
Thus, if we reorder the optimal interpolation 
matrix so that it has the form
\begin{equation} \label{eqn-P-reordered}
\Popt= \left[
    \begin{array}{c}
      P_f \\
      P_c
    \end{array}
    \right], 
\end{equation}
and such that $P_{c}$ is non-singular, then it follows that
the interpolation matrix
\begin{equation} \label{eqn-P-postscaled}
\PoptP = \Popt P_c^{-1} = \left[
    \begin{array}{c}
      \bar{W} \\
      I
    \end{array}
    \right] , \quad  \bar{W} = P_f P_c^{-1}
\end{equation}
also minimizes $\sup\limits_{\Xv\neq 0}\kappa(\Xv)$. 
\end{remark}

The optimal $P$ given in Lemma~\ref{lem:optimalP} and the resulting classical AMG form, $\PoptP$ are in general 
not of direct use in practice since they require the computation of $n_{c}$ eigenvectors of the generalized eigenproblem, given by ~\eqref{eq:gep}, and they yield a dense interpolation matrix. In the remainder of this section, we make various  
 connections between the optimal classical AMG interpolation matrix 
and existing two-grid theory used in deriving classical AMG forms of interpolation.

\subsection{An approximation property and ideal interpolation}
AMG approaches for constructing interpolation $P$ are based on 
an approximation property of the coarse space, 
which is formulated as
\begin{equation}\label{eq:wap}
    \mu_{X}(PR) := \operatorname{sup}_\Xv\dfrac{\|(I - PR)\Xv\|_X^2}{\|\Xv\|^2_A} \leq \eta \ \forall \Xv \in \mathbb{R}^{n},
\end{equation}
where $R: \mathbb{R}^{n} \rightarrow \mathbb{R}^{n_{c}}$ defines the coarse variable type, i.e,. $u_c = Ru$, and must be chosen such that 
$RP = I$ and the matrix $X\in \mathbb{R}^{n\times n}$ is symmetric positive
definite. 

\begin{remark}
Note that the left side 
in~\eqref{eq:wap} will precisely
determine the convergence rate if $X=\XMt$ and $R=(P^T\XMt P)^{-1} P^T \XMt$.
If $X$ is not equal to $\XMt$, but instead bounded from above such that
\begin{equation}\label{eqn:bnd}
(\XMt\Xv,\Xv) = \|\Xv\|^2_*\le \sigma(X\Xv,\Xv),
\end{equation}
then 
\begin{eqnarray*}
&&\|(I-\Pi_\XMt)\Xv\|_\XMt^2 \le \|(I-PR)\Xv\|_\XMt^2 \le \sigma\|(I-PR)\Xv\|_X^2 .
\end{eqnarray*}
\noindent
As a consequence, the two-level method is a uniform contraction in
$\|\cdot\|_A$-norm \emph{if} $\eta$ is uniformly bounded 
for some $X$ such that \eqref{eqn:bnd} holds.  
 A typical choice for $X$, which motivates the classical
AMG approach, is $X=D$ (the diagonal of $A$). 
\end{remark}

As noted above, in the classical AMG setting that is the focus of this paper, interpolation has the form given in
\eqref{eq:classP}.
The choice of restriction $R$ used in defining the coarse variable type in this setting is not unique in that it only needs to satisfy
$\ZP P =I$.  For example, an obvious choice (and the one used most often in practice)
 is given by {\em injection}, that is, by setting 
\begin{equation}
\ZP = \left[
    \begin{array}{cc}
      0 & I
    \end{array}
    \right] \quad \text{such that} \quad u_c = R u = u|_{C},
    \end{equation}
    implying that $u_c$ is the restriction of $u$ to the set of indices in $C$.  
 Though this leads to practical measures for estimating the quality of $C$, it in 
 general gives only an upper bound on $\kappa_\sharp$ even for $X=\XMt$.
 On the other hand, if we set $R=(P^T\XMt P)^{-1} P^T \XMt$, then we precisely
 recover the optimal constant $\kappa_\sharp$ in~\eqref{eq:sharp}, but this choice of $R$ is more 
 difficult to handle in practice.  
 Moreover, this choice will not give the classical AMG form of optimal interpolation, since minimization over 
 $P$ will give $\Popt$ as the minimizer (by Lemma 1), which does not have the classical AMG form of $P$ given in~\eqref{eq:classP}.   In particular, for this choice of $R$, the product $R \Popt = I$, but $R \PoptP \neq I$.
 For the optimal classical AMG case, 
 the proper choice of  $R$ is given by $\RoptR$ defined in Lemma 3 below.  
 The derivation of the result uses the following generalization of a lemma from~\cite{Falgout} that gives the optimal choice for $P$ given some $R$ such that $\ZP P =I$ and any SPD matrix $X$.  The choice here is more general in the use of $\Zid$ which is any full rank matrix of dimension $n_c$, where for the specific choice $\Zid= R^T$ we arrive at the formulation given in~\cite{Falgout}.

\begin{lemma}\label{lem:optPRS}
    Let $R$ be given and satisfy $\ZP P =I$.  Define $S: \mathbb{R}^{n_{f}} \rightarrow \mathbb{R}^{n}$, where $n_{f} = n - n_{c}$ 
    such that $RS = 0$. Then, ideal interpolation is the solution to the $\min$-$\max$ problem $\displaystyle{ \XMx_X= \min_{P} \mu_{X}(PR)}$
   and is given by
    \begin{equation}\label{eq:idealP}\displaystyle{
    \Pid: =       (I-\Pi_A(\ZS))Z_{\ast} = 
(I - \ZS (\ZS^T A \ZS)^{-1} \ZS^T A) \Zid} 
    \end{equation}
    where $\Pid$ satisfies (by definition) $\Pid^{T}AS = 0$ and $\Zid: \mathbb{R}^{n_{c}} \rightarrow \mathbb{R}^{n}$ is full rank. Further, the corresponding minimum has the closed form
    \begin{equation}\label{eq:minOP}
   \XMx_X = \frac{ 1 }{ \lambda_{\min} ( (\ZS^T X \ZS)^{-1} (\ZS^T A \ZS) ) }. 
    \end{equation}
\end{lemma}

In order to minimize the measure $ \mu_{X}$
with respect to $P$, the only requirement that ideal interpolation $\Pid$ needs to satisfy, apart from the usual requirements that $\ZP \Pid=I$ and $\ZP \ZS=0$, is the condition that $\Pid^TA\ZS=0$ (see Equation (3.9) in \cite{Falgout}).  This condition in turn is satisfied for the choice above involving the $A$-orthogonal projection on $\operatorname{range}(S)$.  Moreover, if we consider the condition $\ZP \Pid=I$, we find that 
\begin{equation}\label{eq:RP-RZ}
\ZP \Pid= \ZP(I-\Pi_A(\ZS))Z_{\star}=\ZP Z_{\star} ,
\end{equation}
(since $\ZP \ZS=0$) and, thus, assuming that $\ZP \Pid = I$ is equivalent to the assumption that $\ZP Z_{\star}=I$.  Now, if we choose $Z_{\star}=\ZP^T$, then it must hold that $\ZP \ZP^T=I$, which is the result that was used in~\cite{Falgout} for the construction of $\Pid$ with $Z_{\star} = \ZP^T$.  However, in our construction we consider
the definition of $\Pid$ for $\Zid$ which allows for more general choices of $R$, as well as the relation between $\Pid$ and $\Popt$ that we derive next, which does not hold
when we assume $\Zid = R^T$. 

The next lemma, Lemma~\ref{lem:sharpGAMG}, concerns the choice of $Z_{\star}=\ZP^T$ and shows that our format is indeed more general in terms of the
choice of $\Zid$. 
Specifically, we show that with this more general choice of $\Zid$ (one that is not directly related to the choice of $R$) the ideal and optimal forms of interpolation are one in the same, given proper choices of $R$, $S$ and $Z_{\star}$.
%
%
%
\noindent
The minimizer in this lemma gives a generalization of the so-called ideal interpolation matrix arising in classical AMG.
Next, we derive a result relating this form of ideal interpolation to the optimal one.

\begin{remark}
For the upper ideal $\Pid$ we can also see it can also has following format:
\begin{equation}
\Pid = \left[ \begin{array}{cc} \ZS & \Zid \end{array} \right] \left[
    \begin{array}{c}
      -(\ZS^TA\ZS)^{-1}(\ZS^TA\Zid) \\
      I
    \end{array}
    \right]. 
\end{equation}
which clearly show that when your choice of $\left[ \begin{array}{cc} \ZS & \Zid \end{array} \right]= I$, then the ideal $\Pid$ will have classical AMG format~\eqref{eq:classP}. In ~\cite{Falgout} there is a specific case when choosing $\Zid = R^T$.
\end{remark}

\begin{lemma}\label{lem:sharpGAMG}
    For the optimal classical AMG form of interpolation $P= \PoptP$ and for $X = \XMt$ one obtains, in light of Lemmas~\ref{lem:optimalP} and~\ref{lem:optPRS}, that 
    \begin{equation}\label{eq:optCR}
    \RoptR =(\PoptP P_c P_c^T)^T \XMt, \quad \Sopt = \begin{pmatrix} \Xv_{n_{c}+1} & \cdots & \Xv_{n}\end{pmatrix} 
    \end{equation}
    and
    \begin{equation}
     \XM^{\sharp}_{\XMt} := \frac{ 1 }{ \lambda_{\min} ( (\Sopt^T \XMt \Sopt)^{-1} (\Sopt^T A \Sopt) ) }=\frac{1}{\lambda_{n_c+1}}. \end{equation}
For these choices of $R= \RoptR$, $S=\Sopt$ and $P=\PoptP$, we have that $\RoptR \PoptP = I$, $\RoptR \Sopt = 0$ and $ \PoptP \RoptR = \PiMt(\Popt) = \Pi_{A} (\Popt),$ so that the assumptions
of Lemma~\ref{lem:optPRS} hold.  
Moreover, the resulting ideal P obtained from minimizing the measure in \ref{eq:wap} has the same form as the optimal one, namely
\begin{equation}\label{eq:idVop} 
\Pid = (I - \Sopt (\Sopt^T A \Sopt)^{-1} \Sopt^T A) Z_{\star}=\PiMt(\PoptP)Z_{\star}=\PoptP (\RoptR Z_{\star})=\PoptP,
\end{equation}
where $\Pid$ is defined in Lemma \ref{lem:optPRS} and from \eqref{eq:RP-RZ} we require $\RoptR Z_{\star} = I$.  We postpone further discussion on the choice of $\Zid$ to the remark below that follows the proof of this lemma.  

\begin{proof}
The result in \eqref{eq:idVop} is established by noting that $( \Xv_{i},\Xv_{j})_{\XMt} = \delta_{ij}$  and, thus, $\Sopt^{T}\XMt \Sopt = I$ 
and
$$\Sopt^{T}A\Sopt=\operatorname{diag}(\lambda_{n_{c}+1},\ldots,\lambda_{n})=\Lambda_{\Sopt},$$
implying
  \[
        \Sopt^T\XMt= (\XMt \Sopt)^T = (A \Sopt \Lambda_{\Sopt}^{-1} )^T= \left(\Sopt^{T}A\Sopt\right)^{-1}\Sopt^TA.
 \]
Hence, 
\[
\PiMt(\Sopt)=\Sopt (\Sopt^T \XMt \Sopt)^{-1} \Sopt^T\XMt = \Sopt \Sopt^T \XMt=\Sopt \left(\Sopt^{T}A\Sopt\right)^{-1}\Sopt^TA=\Pi_{A}(\Sopt).
\]
Likewise, one can show that $\PoptP\RoptR=\PiMt(\PoptP)=\Pi_{A}(\PoptP)$.
 Now, since $I-\PiMt(\Sopt)=\PiMt(\PoptP)$, which follows by definition, we have
\[
\Pid = (I - \Sopt (\Sopt^T A \Sopt)^{-1} \Sopt^T A) \Zid=\PiMt(\PoptP)\Zid=\PoptP (\RoptR \Zid) = \PoptP,
\]
where $\RoptR \Zid = I$ is by equivalence ~\eqref{eq:RP-RZ}. 
Thus, if one chooses the optimal forms $R = \RoptR$, $S = \Sopt$ and $X = \XMt$, then the measure in \eqref{eq:wap} also has as its minimizer the optimal form of interpolation $\Popt$ given in Lemma 1. 
\end{proof}
\end{lemma}

\begin{remark}
For the standard case where $\Zid=\RoptR^T$, it follows that $\RoptR \RoptR^T \neq I$ and so $\RoptR \Pid \neq I$ and the results from the previous lemma do not hold.
Indeed, based on definition of $\RoptR$ we have 
$$  \RoptR \Pid: = \RoptR  (I-\Pi_A(\ZS))\RoptR^T= \RoptR \RoptR^T=(\PoptP P_c P_c^T)^T \XMt \XMt (\PoptP P_c P_c^T)= P_c\Popt^T\XMt^2 \Popt P_c^T \neq I.$$
Hence, from this example we see that the $\Pid$ construction in~\cite{Falgout} is not general enough for the derivation in the previous lemma since it does not apply to this optimal choice of $\ZP=\RoptR.$
\end{remark}

The next result characterizes the standard so-called ideal interpolation matrix arising in the classical AMG literature.
\begin{remark}                           
Given the classical AMG form of matrix-dependent interpolation in
~\eqref{eq:classP} and the splitting $C \subset  \{1,...,n\}$ and 
$F =  \{1,...,n\} \setminus C$, consider reordering the matrix such that
\begin{equation}\label{eq:camgA}
A = \left[\begin{matrix}
    A_\ff & A_\fc  \\
   A_\cf & A_\cc 
\end{matrix} \right]
\begin{array}{lr}
\} \: F \\ \} \:  C
\end{array} .
\quad
\end{equation}
Then, practical choices of $R$ and $S$ are as follows
\begin{equation}\label{eq:rns}
\ZP = \left[
    \begin{array}{cc}
      0 & I
    \end{array}
    \right] 
\quad \text{and} \quad
\ZS = \left[
    \begin{array}{c}
      I \\
      0
    \end{array}
    \right] ,
\end{equation}
where for these choices $RS=0$ and $RP = I$.
In this setting, we have that ideal interpolation is given by 
\begin{equation}\label{eq:camgP} 
\Pid = 
 \left[\begin{matrix}
   W_\star \\ I 
\end{matrix} \right] , \quad W_\star = -A_\ff^{-1} A_\fc \: .
\end{equation}
\end{remark}
Hence, for these simpler choices of $R$ and $S$, the optimal and ideal $P$ are in general different and, as we show
in the numerics section, this difference can lead to substantial changes in convergence rates of the resulting
two-grid method in certain cases.   We note though that this form of 
the ideal $P$ is also optimal in terms of energetic stability in the following sense.  

\begin{remark}
In ~\cite{Falgout} the following is derived: Assuming that the coarse variables have been constructed so that $\XMx_{\XMt}$ is bounded for all
$\Xv \neq 0$, then using a $P$ that satisfies
the following stability property also implies convergence 
of the resulting two-level method
\begin{eqnarray}
& &
\label{eqn-harmonic-P}
\XIProd{A P\ZP \Xv}{P\ZP\Xv} \le \beta \XIProd{A\Xv}{\Xv}
\quad \forall \: \Xv,
\rule{0em}{3ex}
\end{eqnarray}
where $\beta \ge 1$ is a constant. 
This more general result is interesting because it  
allows for various approaches of defining interpolation.
Moreover, it separates the tasks of selecting the coarse variables
and defining interpolation.  We note that \eqref{eqn-harmonic-P}
has been used extensively in the literature~\cite{Mandel98energyoptimization,Wan,Brannick_Trace_06} to derive various
techniques for constructing an energetically stable $P$.  
It is also easy to see that if we choose the optimal forms of $R$ and $S$ given in Lemma~\ref{lem:sharpGAMG}, that is, 
$\RoptR =(\PoptP P_c P_c^T)^T \XMt$ and $\Sopt = \begin{pmatrix} \Xv_{n_{c}+1} & \cdots & \Xv_{n}\end{pmatrix}$ such that
$\XMx_{\XMt} = \frac{1}{\lambda_{n_c+1}}$, then with the optimal form of classical AMG interpolation $\PoptP$ given in~\eqref{eqn-P-postscaled}
it follows that $\PoptP \RoptR = \Pi_{\XMt} (\PoptP) = \Pi_A (\PoptP)$ (see the proof in the lemma) so that $\beta = \beta_\sharp = 1$ in~\eqref{eqn-harmonic-P}, which is of course the minimal value since $\Pi_A (\PoptP) := \PoptP (\PoptP^TA\PoptP)^{-1}\PoptP^TA$ is the $A$-orthogonal projection on $\operatorname{range}(\PoptP)$.  
\end{remark}

\section{Compatible Relaxation and Optimal Interpolation}\label{sec:sharpCR}
In this section, we introduce the notion of compatible relaxation (CR) and show how the process can be used 
to obtain a bound on the convergence rate of the two-grid method with ideal interpolation, namely
we present a bound on $\XMx_{\XMt}$ that uses the convergence rate of CR.
The resulting bound in turn gives a measure of the quality of the coarse variable set, since
it shows that, for the given set $C$ there exists an interpolation matrix (the ideal one) such that the 
two-grid method convergence is acceptable.  We then consider CR in terms of $\RoptR$ and
$\Sopt$ defined in Lemma~\ref{lem:sharpGAMG} and contrast this approach to the more practical choices that have been considered 
in our previous studies of this topic.

Compatible relaxation, as defined by Brandt \cite{PBrandt_2000}, is
{\em a modified relaxation scheme that keeps the coarse-level
variables invariant}.  
Generally, the
compatible relaxation iteration is defined as 
\begin{equation} \label{eqn-comp-relax-SF}
    \Xv_{k+1} = (I - S M_S^{-1} S^TA ) \Xv_{k}, 
\end{equation}
where $M_S = \ZS^T M \ZS$ and, as in Lemma~\ref{lem:optPRS}, we assume that $R$ and $S$ are chosen such that
$RS=0$.  Note that from this assumption it follows that
$R \Xv_{k+1} = R \Xv_k$ 
and thus we can consider compatible relaxation in the $L_2$ complementary space, leading to an iteration 
of the form
\begin{equation} \label{eqn-comp-relax-S}
    \Xv_{k+1} = (I -  M_S^{-1} A_S ) \Xv_{k}, 
\end{equation}
where $A_S = S^TAS$.  Now, 
the convergence rate of this iteration 
is related to the measure $\XMx_{\XMt}$ in \eqref{eq:minOP}
as follows
\begin{equation} \label{eqn-cr-convergence}
\XMx_{\XMt} \leq\frac{\Delta^2}{1 - \rho_S} .
\end{equation}
Here, $\Delta \geq 1$ measures the deviation of $M$ from its symmetric
part (see~\cite{Falgout}) and 
\begin{equation} \label{eqn-rho-s}
\rho_S = \XNorm{(I - M_S^{-1} A_S)}_{A_S} .
\end{equation}  Note that, although we
use $\rho$ to represent the spectral radius of a matrix, the quantity
$\rho_S$ is in general only an upper bound for the spectral radius of
compatible relaxation; it is equal to the spectral radius when $M$ is
symmetric.  For this reason, we work with the symmetrized Gauss Seidel smoother
$\XMt$ in all cases, which gives $\Delta = 1$
in the above bound.  

If iteration \eqref{eqn-comp-relax-S} is fast to converge, then 
$\XMx_{\XMt}$ is bounded, that is, fast convergence of CR implies a coarse variable set of good quality 
and the existence of
a $P$ such that the resulting two-level method is uniformly convergent.  
One can then estimate the value of $\rho_s$ in \eqref{eqn-rho-s} in practice by running 
the compatible relaxation iteration in \eqref{eqn-comp-relax-S} and monitoring its convergence.  
\begin{remark}
In the classical AMG setting, given 
a choice of the coarse and fine variable sets, $C$ and $F$, and assuming that $R$ and $S$ are defined
as in~\eqref{eq:rns} we have that
the iteration in \eqref{eqn-comp-relax-S} 
reduces to simple $F$-relaxation
\begin{equation} \label{eqn-Frelax}
\Xv_{k+1} = (I - M_\ff^{-1} A_\ff) \Xv_{k}, 
\end{equation}
which is straightforward to compute.
We note that though this variant of CR is user-friendly, it has been observed in practice
that its spectral radius does not provide an accurate prediction of the convergence
rate of the two-level method with ideal interpolation in some cases~\cite{Brannick_Trace_06}.  
Thus, the bound in~\eqref{eqn-cr-convergence} may not be sharp for such problems.
\end{remark}

\begin{remark}
When $P = \Pid$ (the ideal interpolation operator) and $R$ and $S$ are defined as in \eqref{eq:rns}, then it is easy to show that 
\begin{equation}\label{eq:SharpCR}
I - \Pi_A(\Pid) = I-\Pid(\Pid^{T}A\Pid)^{-1}\Pid^{T}A=S(S^{T}AS)^{-1}S^TA= \begin{pmatrix}
I & -W\\
0 &  0
\end{pmatrix},
\end{equation}
where $W=W_\star$ is given in~\eqref{eq:camgP}.  
This result implies that the spectral radius of the two-grid method with ideal interpolation, $\rho(E_{TG}(\Pid)) $ with
$E_{TG}(\Pid) := (I-M^{-1}A)(I-\Pi_A(\Pid))$, can be accurately estimated in practice if an estimate of $A_\ff^{-1}$ 
is available.  Moreover, the fast convergence of the $F$-relaxation form of CR  in  \eqref{eqn-Frelax} implies
that $A_\ff$ is well conditioned and (as we show later in the numerical experiments) that $A_\ff^{-1}$ can be efficiently estimated using 
a polynomial approximation~\cite{Brannick_Trace_06}.  
In particular, if we consider compatible relaxation defined by \eqref{eqn-comp-relax-S}
and assume that   
$\rho_s \le \delta < 1, 
$
then we have that 
\begin{equation}\label{eq:condest}
\kappa(A_\ff) \le \kappa(M_\ff)\frac{1+\delta}{1-\delta},
\end{equation}
where $\kappa(A)$ is the condition number of the matrix $A$.
In the numerical tests, we use the diagonally scaled preconditioned Conjugate 
Gradient iteration to compute an approximation to $A_\ff^{-1}$.
\end{remark}
\noindent
The approximation obtained by CR can be made optimal if we use the optimal forms of $R$, $S$, and $P$, 
as defined in Lemma 3.  
\begin{remark}  If, instead, we use the optimal choices of $R$ and $S$,  namely, $\RoptR$ and $\Sopt$, and we assume
that the smoother is symmetric Gauss Seidel, then
instead of arriving at the $F$-relaxation form of CR given in \eqref{eqn-Frelax} we have that
\begin{equation}\label{eqn:CRopt}
\XMt_{\Sopt} =\Sopt^T\XMt\Sopt = I \quad \text{and} \quad A_{\Sopt} =\Sopt^T A\Sopt = \Lambda_{\Sopt} := \text{diag}(\lambda_{n_c+1}, ..., \lambda_n), 
\end{equation}
where $\lambda_i$ are the generalized eigenvalues of the eigenproblem involving $(A,\XMt)$.  
Thus, the spectral radius of the CR error propagation matrix is given by
\begin{equation}\label{eqn:CRoptRho}
\rho_{\Sopt}(I - \XMt_{\Sopt}^{-1} A_{\Sopt}) = 1-  \min_i (\Lambda_{\Sopt})_{ii} = 1 - \lambda_{n_c+1}.
\end{equation}
Here, unlike the simplified version of CR given above, the convergence rate of this form of CR does not depend on the actual
coarse points that are chosen, but only on the cardinality of the coarse variable set $n_c = |C|$.
However, if we recall the optimal form of classical AMG interpolation given in \eqref{eqn-P-postscaled}:
\begin{equation} 
\PoptP = \Popt P_c^{-1} = \left[
    \begin{array}{c}
      \bar{W} \\
      I
    \end{array}
    \right] ,  \quad \bar{W} = P_f P_c^{-1},
\end{equation}
then we observe that $P_c$ used in constructing the classical AMG form of optimal interpolation
does depend on the choice of 
coarse grid points.  Hence, in this setting we can use CR to determine the number of coarse points 
that are required in order to achieve a certain convergence rate of the resulting two-grid method and
then we choose the set $C$ so that $P_c$ is well conditioned.  This assumption in turn implies
$P_c^{-1}$ is computable in the approximation to $\PoptP$ and as we show numerically this also 
tends to produce a local (decaying) optimal interpolation matrix.  

\end{remark}


\subsection{An equivalence between ideal and optimal interpolation for $F$-relaxation}
If we consider the reduction-based AMG framework and assume the forms of $R$ and $S$ given in~\eqref{eq:rns}, then we can establish an equivalence between the minimizers of $\XMx_{\XMt}$ in~\eqref{eq:minOP} and $\kappa_\sharp$ in~\eqref{eq:sharp} with respect to $P$.
However, as we show in the numerical results section, in the case of full smoothing (i.e., using $M$ not $M_\ff$ in the solve phase)
$\XMx_{\XMt}$  does not generally provide a useful estimate of  $\kappa_\sharp$.

The reduction-based AMG method (or algebraic hierarchical basis scheme) is based on the 
following block factorization 
of $A$ in terms of $R$ and $S$:
\begin{equation}
A = \left[\begin{matrix}
   A_\ff & A_\fc  \\
   A_\cf & A_\cc 
\end{matrix} \right] \begin{array}{lr}
\} \: F \\ \} \:  C
\end{array} =
\left[\begin{matrix}
  I & 0 \\
    A_\cf A_\ff ^{-1}   & I 
\end{matrix} \right] 
\left[\begin{matrix}
   A_\ff & 0  \\
  0 & A_\cc - A_\cf A_\ff ^{-1}A_\fc
\end{matrix} \right] 
\left[\begin{matrix}
   I & A_\ff ^{-1}A_\fc   \\
   0 & I
\end{matrix} \right]
, 
\end{equation}
implying
\begin{equation}
A^{-1} =
\left[\begin{matrix}
   I & -A_\ff ^{-1}A_\fc   \\
   0 & I
\end{matrix} \right]
\left[\begin{matrix}
   A_\ff^{-1} & 0  \\
  0 & S_A^{-1}
\end{matrix} \right] 
\left[\begin{matrix}
  I & 0 \\
    -A_\cf A_\ff ^{-1}   & I 
\end{matrix} \right] , 
\end{equation}
where $ S_A = A_\cc - A_\cf A_\ff ^{-1}A_\fc $.
Equivalently, one can solve the system using the two-level method with ideal interpolation $\Pid$, which gives 
$A_c = \Pid^TA\Pid = S_A$ and again results in an exact solve, assuming exact $F$-relaxation, i.e., $M_\ff = A_\ff$ in
\begin{equation}\label{eq:idP}
E_{TG}  = (I - S M^{-1}_\ff S^T A) ( I - \Pi_{A}(\Pid))(I -S M_\ff^{-T} S^T A) = I - B_{TG}^{-1}A,
\end{equation}
where $M_\ff = S^TMS$ and $\Pi_{A}(\Pid)$ is an $A$-orthogonal projection on $\operatorname{range}(\Pid)$
as defined in~\eqref{eq:camgP}.

Now, consider the block factorized smoother
\begin{equation}
M_{HB} = 
\left[\begin{matrix}
  M_\ff & 0 \\
 A_\cf  & \tau I 
\end{matrix} \right]
\left[\begin{matrix}
 I & M_\ff^{-1} A_\fc \\
  0  & I
\end{matrix} \right], 
\end{equation}
where for sufficiently large $\tau$
we have that $M_{HB}^T + M_{HB} - A$ 
is symmetric and positive definite.
We note that 
\begin{equation}
M_{HB}^{-1} = 
\left[\begin{matrix}
  M_\ff & A_\fc \\
  A_\cf  & \tau I + A_\cf M_\ff^{-1} A_\fc
\end{matrix} \right]^{-1} =  
\left[\begin{matrix}
M_\ff^{-1}+\tau^{-1} M_\ff^{-1} A_\fc A_\cf M_\ff^{-1} & -\tau^{-1}M_\ff^{-1}A_\fc \\
-\tau^{-1} A_\cf M_\ff^{-1} & \tau^{-1} I 
\end{matrix} \right] 
\rightarrow
\left[\begin{matrix}
M_\ff^{-1} & 0 \\
0 & 0 
\end{matrix} \right],
\end{equation}
where in the last step we take $\tau \rightarrow \infty$.  Thus, as $\tau$ approaches 
infinity the factorized smoother converges to simple $F$-relaxation involving $M_\ff^{-1}$.  

Next, we set $M_\ff = A_\ff$ and note that in this case the symmetrized smoother is given by:
\begin{equation}
\widetilde{M}_{HB} = 
\left[\begin{matrix}
A_\ff & A_\fc \\
A_\cf & \tau^2 X^{-1} + A_\cf A_\ff^{-1} A_\fc 
\end{matrix} \right] , 
\end{equation}
where $X^{-1} = 2 \tau I - S_A$.
Now, considering the generalized eigen-problem
involving $A$ and $\widetilde{M}_{HB}$, it follows 
that 
\begin{equation}\label{eq:MHBinvA}
\widetilde{M}_{HB}^{-1} A = 
\left[\begin{matrix}
A_\ff & A_\fc \\
A_\cf & \tau^2 X^{-1} + A_\cf A_\ff^{-1} A_\fc 
\end{matrix} \right]^{-1}
 \left[\begin{matrix}
   A_\ff  & A_\fc \\
  A_\cf & A_\cc
\end{matrix} \right] 
=
\left[\begin{matrix}
I & \ast \\
0 & \ast 
\end{matrix} \right] , 
\end{equation}
so that $\lambda_{n_c+1} = \hdots = \lambda_n =1$
for any choice of hierarchical splitting.  
Thus, if we consider using the optimal interpolation 
operator $\PoptP$ given in~\eqref{eqn-P-postscaled},then the above result shows that 
the resulting two-level method  with exact $F$-relaxation (where $M_\ff=A_\ff$ in $\widetilde{M}_{HB}$) gives an exact solve for any value of $\tau$,
since in this case
\[
\|E_{TG}\|_A^2 = 1 - \lambda_{n_{c}+1} = 1-1=0 \quad \text{for any value of $n_c$}
.\]
Also ~\eqref{eq:MHBinvA} also implies $\Sopt=\left[\begin{matrix}
I  \\ 0  
\end{matrix} \right] $ as in ~\eqref{eq:rns}. Then
\begin{equation}
(I-\widetilde{M}_{HB}^{-1} A)^TA\Sopt =  (\widetilde{M}_{HB}-A)\widetilde{M}_{HB}^{-1} A\Sopt=
\left[\begin{matrix}
0 & 0 \\
0 & \ast
\end{matrix} \right]
 \left[\begin{matrix}
  I & \ast \\
  0& \ast
\end{matrix} \right] 
\left[\begin{matrix}
I  \\ 0  
\end{matrix} \right] =0.
\end{equation}
Thus $ \operatorname{range}(I-\widetilde{M}_{HB}^{-1} A)= \operatorname{range}(\Popt)$, since $\Popt$ and $\Sopt$ are $A$-orthogonal.
Then $\|E_{TG}\|_A=0.$


Moreover, since when $\tau \rightarrow \infty$, it follows that $M_{HB}$ 
converges to $F$-relaxation and assuming that $M_\ff = A_\ff$, it follows that  
$$\lim_{\tau \rightarrow \infty} (\PoptP)_\tau \rightarrow \Pid ,$$
that is the optimal interpolation operator converges to the ideal 
one.  Here, we emphasize the dependence of $\PoptP$ 
on $\tau$, since the eigenvectors of the generalized eigen-problem 
involving $\widetilde{M}_{HB}$ depend on $\tau$ and, thus, so do 
the columns of $\PoptP$.
However, when $M_\ff \neq A_\ff$, but instead $M_\ff \approx A_\ff$, then
$\lim_{\tau \rightarrow \infty}  \bar{P}_\tau \neq  \Pid$.  Similarly, as we show in 
the numerical results in the next section, in the case of full smoothing the two-level
method that uses optimal interpolation will generally converge faster than the 
two-level method with ideal interpolation.

For the limiting case of $\tau\rightarrow \infty$, 
we have
 \begin{equation}
 \lim_{\tau\to\infty}M_{HB}^{-1} =\left[\begin{matrix}
M_\ff^{-1} & 0 \\
0 & 0 
\end{matrix} \right],
\end{equation}
We note that for this limiting case, the matrix becomes singular and so a pseudo inverse is needed. Now, 
if $M_\ff^{-1} =A_\ff^{-1} $, then
\begin{equation}
\lim_{\tau\to\infty}\widetilde{M}_{HB}^{\dagger} A = 
\left[\begin{matrix}
I & A_\ff^{-1}A_\fc \\
0 & 0 
\end{matrix} \right] , 
\end{equation}
So in analogy to the case $\tau<\infty$, the optimal P is composed of $n_c$ eigenvectors of the above matrix on the right-hand side, namely, the eigenvectors corresponding to zero eigenvalues of
\begin{equation}
\left[\begin{matrix}
I & A_\ff^{-1}A_\fc \\
0 & 0 
\end{matrix} \right] 
\left[\begin{matrix}
v_{f} \\ v_{c}\end{matrix} \right]
=
\left[\begin{matrix}
0 \\ 0
\end{matrix} \right] ,
\end{equation} 
and so
\begin{equation}
v_{f}=-A_\ff^{-1}A_\fc v_{c}.
\end{equation}
Thus, the eigenvectors have the form
$$
v=\left[\begin{matrix}
-A_\ff^{-1}A_\fc \\
I 
\end{matrix} \right] v_{c} = \Pid v_{c}
$$
Hence, the range of the ideal and the optimal interpolation matrices are the same if F-relaxation is used for smoothing. 

\section{Numerical experiments}\label{sec:num}
\begin{figure}[]
\begin{center}
\includegraphics[width=.25\textwidth]{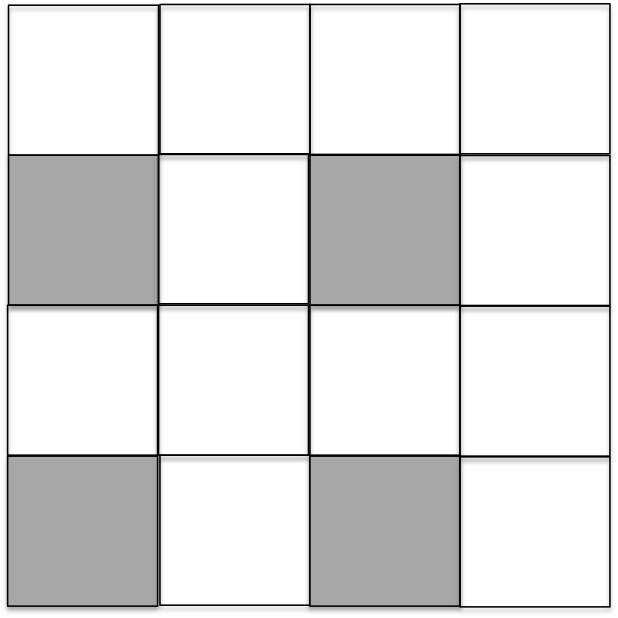}\hspace{.25cm}
\includegraphics[width=.25\textwidth]{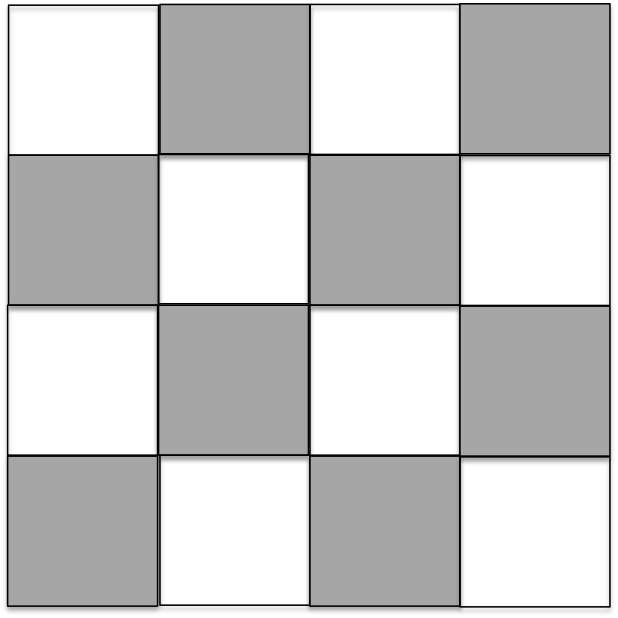}
\caption{Distribution of the jump coefficient $a(x)$. Left: Distribution of P1 and P3; Right: Distribution of P2 and P4.}
\label{fig:jumpdist}
\end{center}
\end{figure}

Numerically, we consider four different test problems coming from discretizations of \eqref{eqn:model}, corresponding to different distributions of the diffusion coefficient $a(x)$.  The first two Problems P1 and P3 are ones where the interfaces of the jumps do not intersect, namely, 
\begin{eqnarray} \label{eqn:P1}
a(x) = 
\begin{cases}
1  \qquad \quad x \in \Omega_1 \; , \\
10^{-k_{ij}}  \quad x \in \Omega \setminus \Omega_1 \; ,
\end{cases}
\end{eqnarray}
where the domain $\Omega_1$ corresponds to the one given by the white regions in the plot on the left in Figure~\ref{fig:jumpdist}.  For problem P1, $k_{ij} = k \in \mathbb{Z}^+$ for all $i,j$, and for problem P3 $k_{ij} \in \{1,2,\hdots,k\}$, where the values are selected randomly with uniform distribution (using built-in Matlab function \texttt{randi}).  In the next two tests, Problems P2 and P4, we consider a checkerboard pattern for the distribution of the jumps, where now $\Omega_1$ corresponds to the white regions in the plot on the right in Figure~\ref{fig:jumpdist}.  For problem P2, the distribution in $\Omega \setminus \Omega_1$ is again uniform with $k_{ij} = k \in \mathbb{Z}^+$ for all $i,j$ and for problem P4 we select the value $k_{ij}$ randomly as in problem P3.  

We use a standard cell-centered finite volume method (see~\cite{eymard2000finite,samarskii2001theory}) for discretizing P1-P4 and choose a structured grid $0 = x_0 < x_1 < \cdots < x_{N+1}$, $x_i = \frac{i}{N+1}$, and $0 = y_0 < y_1 < \cdots < y_{N+1}$, $y_j = \frac{j}{N+1}$.  Note that $h = h_x = h_y = \frac{1}{N+1}$.  Here, each cell $[x_{i-1}, y_{j-1}] \times [x_i, y_j]$ is used as the control volume and the unknowns are located at each cell center $(x_{i- \frac{1}{2}}, y_{j - \frac{1}{2}}) = (x_{i-1} + \frac{h}{2}, y_{j-1}+\frac{h}{2})$.  
To define $a(x)$ on the interfaces of neighboring subdomains we use a harmonic average, e.g., for an interface $S_*$ we have $a^- \neq a^+$ in general, where $a^- = a(x)$ in the volume on one side of the interface $S_*$ and $a^+=a(x)$ in the volume on the other side of $S_*$.  
We then write the discretized system as
\begin{equation}\label{eqn:model-dis}
a_e u_{i+1,j} + a_w u_{i-1,j} + a_n u_{i,j+1} + a_s u_{i,j-1} - (a_e+a_w+a_n+a_s) u_{i,j} = f_{i,j},
\end{equation}
where $f_{i,j} = \int_V f \mathrm{d} V$ and $V$ is the control volume $[x_{i-1}, y_{j-1}] \times [x_i, y_j]$
and $a_*$ are harmonic averages of $a(x)$ on the two neighboring cells as in~\cite{LLM92}. 
We assume Dirichlet boundary conditions and if an edge $S_*$ is on the boundary of $\Omega$, we set $u
_{i+\frac{1}{2}, j} = 0$.  

\subsection{Compatible relaxation: measuring the quality of $C$}
We begin with tests that compare the standard $F$-relaxation form of CR in~\eqref{eqn-Frelax}, the estimate of $E_{TG}(\Pid)$ obtained by using the identity in~\eqref{eq:SharpCR} to apply the coarse-grid correction with ideal interpolation ($P=\Pid$ in~\eqref{eq:camgP}), and the spectral radius of $E_{TG}(\PoptP)$, that is, the two-grid method with optimal interpolation, for Problems P1-P4.  In all tests, we use forward Gauss-Seidel for pre-smoothing and backward Gauss-Seidel for post-smoothing in defining $E_{TG}$ in~\eqref{eq:idP} and we
consider standard full-coarsening ($h \rightarrow 2h$, see Figure~\ref{fig:fc}) to define the coarse variable set $C \subset \{1,...,n\}$ for  different choices of the problem size $n = N\times N$.  The results of these tests for various values of the jump discontinuity defined by parameter $k$ are given in Table~\ref{tab:CR}.  The table at the top contains approximations of the spectral radius of $E_{TG}(\Pid)$ together with the ideal interpolation matrix, i.e., when $ W_\star = 
  -A_\ff^{-1} A_\fc$.  The second table
contains estimates of the two-grid convergence rate obtained by running 5 steps of the iteration~\eqref{eq:SharpCR} starting with a random initial guess, where the action of $A_\ff^{-1}$ in{~\eqref{eq:SharpCR} is approximated by $L=2$ diagonally preconditioned Conjugate Gradient iterations.  
Here, we combine the iteration~\eqref{eq:SharpCR} with Gauss Seidel pre- and post-smoothing, which then mimics the two-grid method
with ideal $P=\Pid$.
The third table contains results of the $F$-relaxation form of CR for symmetric Gauss Seidel, again assuming full
 coarsening as the choice of $C$.  And the bottom table contains results of the two-grid method with optimal interpolation, $\PoptP$.
  
Overall, we see that CR convergence rates are acceptable for all problems and choices of the problem parameters except for Problem P4.  Moreover, from the results obtained in the middle table it is clear that one can obtain an accurate estimate of $E_{TG}(\Pid)$ using~\eqref{eq:SharpCR} with a small number of inner PCG iterations to approximate $A_\ff^{-1}$, again in all cases except for Problem P4.  This poor performance observed for Problem P4 is of course expected since in this case the convergence of $F$-relaxation with full coarsening is $\rho > .9$ and so  $A_\ff$ may not be well conditioned in the sense of the bound given in~\eqref{eq:condest}.
In general, these results suggest that full coarsening may not be a good choice for Problem P4.  However, as we show in the next set of experiments, 
with optimal interpolation full coarsening does give acceptable results. The rapid convergence observed in these tests, especially for Problem P4, can be explained using the results provided in Figure~\ref{fig:spec}, which contains plots of the spectra of $A$ and $(A,\XMt)$. Here, we see that the eigenvalues of $(A,\XMt)$ vary substantially from the eigenvalues of $A$, e.g. in the right plot for Problem P4, less than $\frac{n}{4}$ of the generalized eigenvalues are different from 1. Thus, for our choice of full coarsening with $n=\frac{n_c}{4}$, we obtain very  fast convergence. 

We note that if we instead consider red-black 
coarsening for problem P4, then $F$-relaxation becomes an exact solve, i.e., the spectral radius of the CR iteration, $\rho_s=0$.  In addition, when using the same 
red-black coarsening and iteration~\eqref{eq:SharpCR}, again with 2 inner PCG iterations, we obtain an estimate of the spectral radius of $\rho =.248$ for Problem P4 with $k=8$ and $h =1/2^5$, where the true spectral radius of the two-level method is $.250$ independent of $k$ and $h$.   
Thus, in practice, one can run $F$-relaxation to choose $C$ until the CR iteration converges quickly (so that $A_\ff$ is well conditioned by~\eqref{eq:condest})
and then compute the sharp estimate of $E_{TG}(\Pid)$ defined in terms of ~\eqref{eq:SharpCR} using PCG iterations to approximate the action of $A_\ff^{-1}$ in ~\eqref{eq:SharpCR} in order to obtain a more accurate estimate of the 
convergence rate of the two-grid method using ideal interpolation.  

Finally, in Section 3, we showed that in general the two-level method with optimal interpolation will converge faster than the method with ideal 
interpolation.  Here, we observe these results numerically for Problems P1-P4 in the bottom set of results.   Overall, we observe that the reported convergence results are consistent with the 
theoretical result that $\rho(E_{TG}(\Popt)) \leq \rho(E_{TG}(\Pid))$ and for Problems P2 and P4 a significant improvement is observed.   We note that for the Poisson problem (i.e., $k=0$), the two-grid method with optimal $P$ and full coarsening has its spectral radius bounded by .14 independent of the problem size $n$.  

\begin{table}[ht!]
\begin{center}
Spectral radius of $E_{TG}(\Pid)$
\vspace{.1cm}

{\scriptsize
\begin{tabular}{@{\extracolsep{-1.mm}} |c|cccc|cccc|cccc|cccc|}
 \hline  
       & \multicolumn{4}{c|}{$k=1$} & \multicolumn{4}{c|}{$k=2$} & \multicolumn{4}{c|}{$k=4$} & \multicolumn{4}{c|}{$k=8$} \\ \hline    
    Size   & P1 & P2 & P3 & P4 & P1 & P2 & P3 & P4 & P1 & P2 & P3 & P4  & P1 & P2 & P3 & P4 \\ \hline
$16^2$  & .259&.255 &.300&0.397   &    .251&.251&.297&.535    &      .250&.250&.298&.577   &     .250&.250&.294&.679   \\             
$32^2$  & .260&.256 &.302&0.445   &    .251&.251&.301&.649    &      .250&.250&.293&.791   &     .250&.250&.285&.887   \\                  
$64^2$  & .261&.256 &.303&0.473   &    .251&.251&.301&.714    &      .250&.250&.294&.879   &     .250&.250&.292&.991   \\             
$128^2$& .261&.256 &.305&0.471 Ê &    .251&.251&.301&.729    &      .250&.250&.298&.924   &     .250&.251&.294&.997    \\            
 \hline 
\end{tabular}
}

\vspace{.1cm}
Approximation of $E_{TG}(\Pid)$ using identity \eqref{eq:SharpCR}
\vspace{.1cm}

{\scriptsize
\begin{tabular}{@{\extracolsep{-1.mm}} |c|cccc|cccc|cccc|cccc|}
 \hline  
       & \multicolumn{4}{c|}{$k=1$} & \multicolumn{4}{c|}{$k=2$} & \multicolumn{4}{c|}{$k=4$} & \multicolumn{4}{c|}{$k=8$} \\ \hline    
    Size   & P1 & P2 & P3 & P4 & P1 & P2 & P3 & P4 & P1 & P2 & P3 & P4  & P1 & P2 & P3 & P4 \\ \hline
$16^2$   &.240&.235&.249&.209   &.233&.231&.244&.210   &.232&.231&.239&.217    &.232&.231&.231&.225 \\ 
$32^2$   &.245&.243&.253&.198   &.241&.241&.250&.204   &.240&.241&.247&.205    &.240&.241&.239&.220\\
$64^2$   &.244&.242&.252&.200   &.239&.239&.250&.205   &.239&.239&.247&.216    &.239&.239&.237&.225\\
$128^2$   &.234&.237&.220&.202   &.240&.240&.231&.206   &.240&.240&.236&.214  &.240&.240&.238&.223\\

 \hline 
\end{tabular}
}

\vspace{.1cm}
Compatible Relaxation iteration~\eqref{eqn-Frelax} with symmetric Gauss Seidel
\vspace{.1cm}

{\scriptsize
\begin{tabular}{@{\extracolsep{-1.mm}} |c|cccc|cccc|cccc|cccc|}
 \hline  	`
       & \multicolumn{4}{c|}{$k=1$} & \multicolumn{4}{c|}{$k=2$} & \multicolumn{4}{c|}{$k=4$} & \multicolumn{4}{c|}{$k=8$} \\ \hline    
    Size   & P1 & P2 & P3 & P4 & P1 & P2 & P3 & P4 & P1 & P2 & P3 & P4  & P1 & P2 & P3 & P4 \\ \hline
$16^2$   &.242&.176&.512&.693   &.075&.052&.493&.839   &.007&.005&.499&.937    &7e-5&5e-5&.500&.999 \\ 
$32^2$   &.243&.177&.524&.786   &.075&.052&.520&.939   &.007&.005&.516&.995    &7e-5&5e-5&.512&1.00\\
$64^2$   &.244&.177&.530&.777   &.075&.052&.527&.927   &.007&.005&.522&.989    &7e-5&5e-5&.515&1.00\\
$128^2$   &.244&.178&.533&.790   &.075&.052&.526&.951   &.007&.005&.523&.998  &7e-5&5e-5&.517&1.00\\

 \hline 
\end{tabular}
}
\vspace{.2cm}

Spectral radius of $E_{TG}(\Popt)$
\vspace{.1cm}

{\scriptsize
\begin{tabular}{@{\extracolsep{-1.mm}} |c|cccc|cccc|cccc|cccc|}
\hline  
        & \multicolumn{4}{c|}{$k=1$} & \multicolumn{4}{c|}{$k=2$} & \multicolumn{4}{c|}{$k=4$} & \multicolumn{4}{c|}{$k=8$} \\ \hline    
    Size &   P1 & P2 & P3 & P4 & P1 & P2 & P3 & P4 & P1 & P2 & P3 & P4  & P1 & P2 & P3 & P4 \\ \hline
$16^2$    &.041&.024&.124 &.132	&.005&.002&.108&.120	&5e-5&3e-5&.102 &.102	&5e-9&2.5e-9&.065&.065\\
$32^2$    &.042&.024&.134&.148	&.005&.002&.131&.140	&5e-5&3e-5&.126 &.128	&5e-9&2.5e-9&.087&.117\\
$64^2$    &.042&.024&.137&.154	&.005&.002&.136&.152	&5e-5&3e-5&.132 &.146 &5e-9&2.5e-9&.124&.151\\
$128^2$  &.042&.024&.140&.160	&.005&.002&.139&.159	&5e-5&3e-5&.136 &.157	&5e-9&2.5e-9&.127&.159\\

 \hline 
\end{tabular}
}
\end{center}
\label{tab:CR}
\vspace{.1cm}
\caption{Spectral radius of two-grid methods with ideal interpolation (top), using CG to approximate $A_\ff^{-1}$ (top-middle), results for compatible relaxation (bottom-middle)  and with optimal interpolation (bottom), applied to Problems P1-P4 (with symmetric GS as smoother).}
\end{table}%


\subsection{Compatible relaxation with optimal interpolation}\label{sec:optPCR}

As discussed in Section \ref{sec:sharpCR}, the convergence
rate of the $F$-relaxation form of compatible relaxation can be used to measure the quality of the coarse variable
set in that it bounds the $\min-\max$ solution $\XMx_{\XMt}$ in~\eqref{eq:wap}. Recall that this
construction assumes the so-called ideal interpolation operator in~\eqref{eq:idealP}, which as we showed above
does not give the best choice of the classical AMG form of interpolation due to the forms of $R$ and $S$ given 
in~\eqref{eq:rns} that are assumed in this setting. 

Our aim in this
section is to study the use of CR together with the optimal forms of $R$ and $S$ given in Lemma~\ref{lem:sharpGAMG}, namely, $\RoptR =(\PoptP P_c P_c^T)^T \XMt$ and $\Sopt = \begin{pmatrix} \Xv_{n_{c}+1} & \cdots & \Xv_{n}\end{pmatrix}$, which leads to the optimal interpolation matrix as the 
minimizer of $\XM_{\XMt}$.  Here, $P_c$ defines the matrix used in deriving the classical AMG form of $\Popt$ given in~\eqref{eqn-P-postscaled}:
\begin{equation*} 
\PoptP = \Popt P_c^{-1} = \left[
    \begin{array}{c}
  P_f P_c^{-1} \\
      I
    \end{array}
    \right] , \quad
\Popt= \left[
    \begin{array}{c}
      P_f \\
      P_c
    \end{array}
    \right] \leftarrow \left[ \begin{matrix} \Xv_{1} & \cdots & \Xv_{n_c}
    \end{matrix}  \right] ,
\end{equation*}
with the set $C$ chosen such that $P_{c}$ is non-singular.   The arrow notation is used to denote that $\Popt$ with columns consisting of the smallest 
$n_c$ eigenvectors of the generalized eigenproblem for $(A,\XMt)$ is reordered according to the $C-F$ splitting of the unknowns.   Note that from Lemma~\ref{lem:optPRS}  we have that the classical AMG form of optimal interpolation $\PoptP$ not only yields the optimal two-grid convergence rate, i.e., it gives 
$\kappa_\sharp$ in~\eqref{eq:sharp}, it also minimizes the approximation property $\XM_{\XMt}$ with respect to $P$.  

Given the above choice of $S=\Sopt$ we have that the spectral radius of CR is given by (see~\eqref{eqn:CRoptRho})
\begin{equation}
\rho_{\Sopt}(I - \XMt_{\Sopt}^{-1} A_{\Sopt}) = 1 - \lambda_{n_c+1}.
\end{equation}
Here, $\XMt$ denotes the symmetric Gauss-Seidel smoother.
Note that, the CR rate gives the same result as the convergence rate of the two-grid method with optimal interpolation, namely, 
the same rate obtained with $\kappa_{\sharp}$.  However, unlike the simplified $F$-relaxation version of CR given in~\eqref{eqn-Frelax}, the convergence rate of this form of CR does not depend on the actual
coarse points that are chosen, instead it only depends on the cardinality of the coarse variable set $n_c = |C|$.
The choice of the coarse variable set now defines the matrix $P_c$ that is used in constructing the classical AMG form of optimal interpolation.
Hence, in this setting we can use CR to determine the number of coarse points 
that are required in order to achieve a certain convergence rate of the resulting two-grid method and
then we choose the set $C$ so that $P_c$ is well conditioned. 

The problem of finding a well-conditioned submatrix $P_c$ can be viewed as finding a subset of $n_{c}$ columns of $P_\sharp^T$ that form a
well conditioned basis of $\mathbb{R}^{n_c}$ or in looser terms a set of $n_{c}$ columns that are as linearly independent as possible. Besides the
condition number of the basis, i.e., a measure for its orthogonality, another natural formalization of this question is the volume of the parallelepiped
spanned by the bases. In the case that all columns of $P_{\sharp}^{T}$ have comparable norm, a maximal volume (determinant) basis has a small condition number and vice
versa. Unfortunately, finding the maximal volume submatrix is an NP-hard problem (cf.~\cite{CivrMagd2013}), but there exists a greedy algorithm 
(cf. Algorithm~\ref{alg:maxvol}) that is able to find a sub-matrix $P_c$ with locally maximal determinant (see~\cite{Knut1985,GoreOselSavoTyrtZama2010} for details). 
Even though there is no guarantee that the submatrix found in this way has a small condition number,
we find in practice that the large eigenvalues tend to stay within the same order. 
In contrast, the nearly zero eigenvalues of $P_c$ move away from the origin significantly when the algorithm
is applied to some random initial choice of coarse grid.  

Figures \ref{fig:2d-spec} - \ref{fig:2d-spec3} contain results of the algorithm applied to our scalar diffusion problem in \eqref{eqn:model}
for different choices of the diffusion tensor $a(x)$, namely, the distribution for our test problems P1-P4 with $k=4$. For all tests we use a discretization of
the problem with $35\times 35$ finite volumes. Each figure contains plots of the
choice of coarse variable set $C$, denoted by the black circles as well as the volume and condition numbers of the depicted choice of $C$. We depict the initial (random) choice 
of the $C$ set and the set determined by the \texttt{maxvol} algorithm.
Further, we show a plot of the convergence rate versus
the cardinality of $n_c = |C|$, and plots of three (randomly chosen) columns of the resulting optimal classical AMG form of
interpolation that we obtain for the given coarse grid and then computing $\PoptP$ by multiplying $\Popt$ by
$P_c^{-1}$.  
We observe that, in all tests, the set of coarse variables properly aligns with the choice of $a(x)$ and the given smoother. For example, 
 in Figure \ref{fig:2d-spec1} we naturally obtain standard full coarsening. For the results in Figure \ref{fig:2d-spec3} 
 for Problem P4 using a red-black ordered block smoother we see that the coarse points largely lie on the boundaries of the subdomains, 
as expected.  In addition, the plots
of the convergence rates are given on the top right and allow one to choose $n_c=|C|$ with a guaranteed lower bound on the convergence rate of the resulting two-grid 
method, namely the one obtained with $\PoptP$.  Finally, we note that the columns of the resulting optimal classical AMG form of interpolation are highly localized for all test problems.

\bigskip

\begin{algorithm}[H]
    \caption{{\texttt{maxvol}} coarse variable selection -- {Input:} $\Popt$ and $n_c$, {Output:} $C$, $F$\label{alg:maxvol}}
		\SetAlgorithmStyle
        Choose initial set $C$ such that $|C| = n_c$, set $X := (\Popt)_{C}$ \;
        Calculate $\PoptP := \pi\Popt X^{-1} = \begin{pmatrix} I \\ W\end{pmatrix}$\tcp*{\ reordered according to $C$-$F$}
        \While{$\operatorname{max}_{i, j} |(\PoptP)_{i,j}| > 1$}{
            $(i^{\prime},j^{\prime}) = \operatorname{argmax}\{|(\PoptP)_{i,j}|\}$\tcp*{\ $i > n_{c}$, i.e., $i \notin C$}
            Set $C = (C\setminus \{j^{\prime}\}) \cup \{i^{\prime}\}$\tcp*{\ swap rows $j^{\prime}$ and $i^{\prime}$}
            Update $\PoptP$\tcp*{rank one updates}
		}
	\end{algorithm}

    \begin{figure}[ht!]
        \begin{center}
            \hfill
            \begin{minipage}[b]{.25\textwidth}
                \includegraphics[width=1\textwidth]{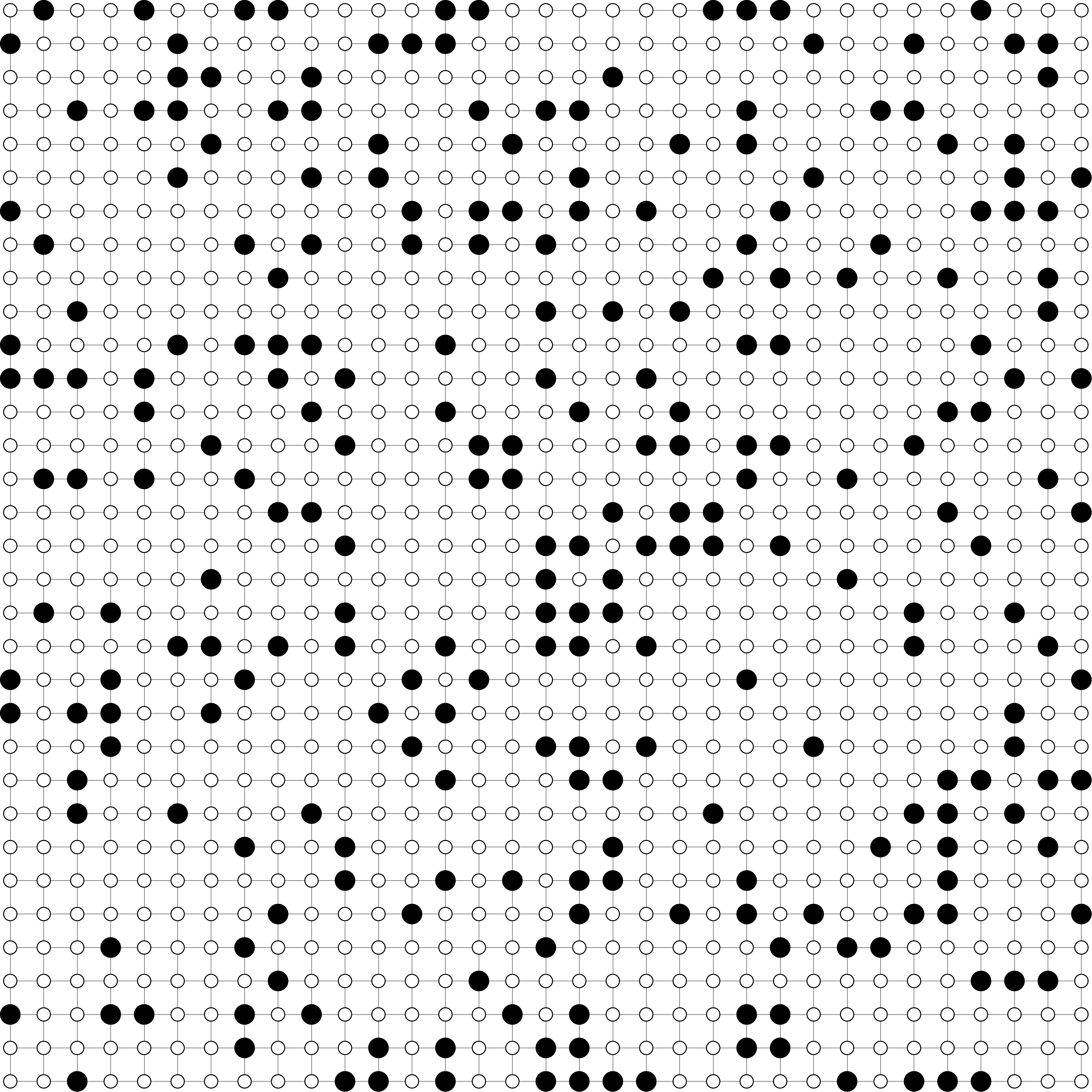}
                {\scriptsize \begin{align*} |\operatorname{det}(P_{c})| & = 2.245\cdot 10^{-245}\\ \kappa(P_{c}) &= 4.709\cdot10^{14}\end{align*}}%
            \end{minipage}
            \hfill 
            \begin{minipage}[b]{.25\textwidth}
                \includegraphics[width=1\textwidth]{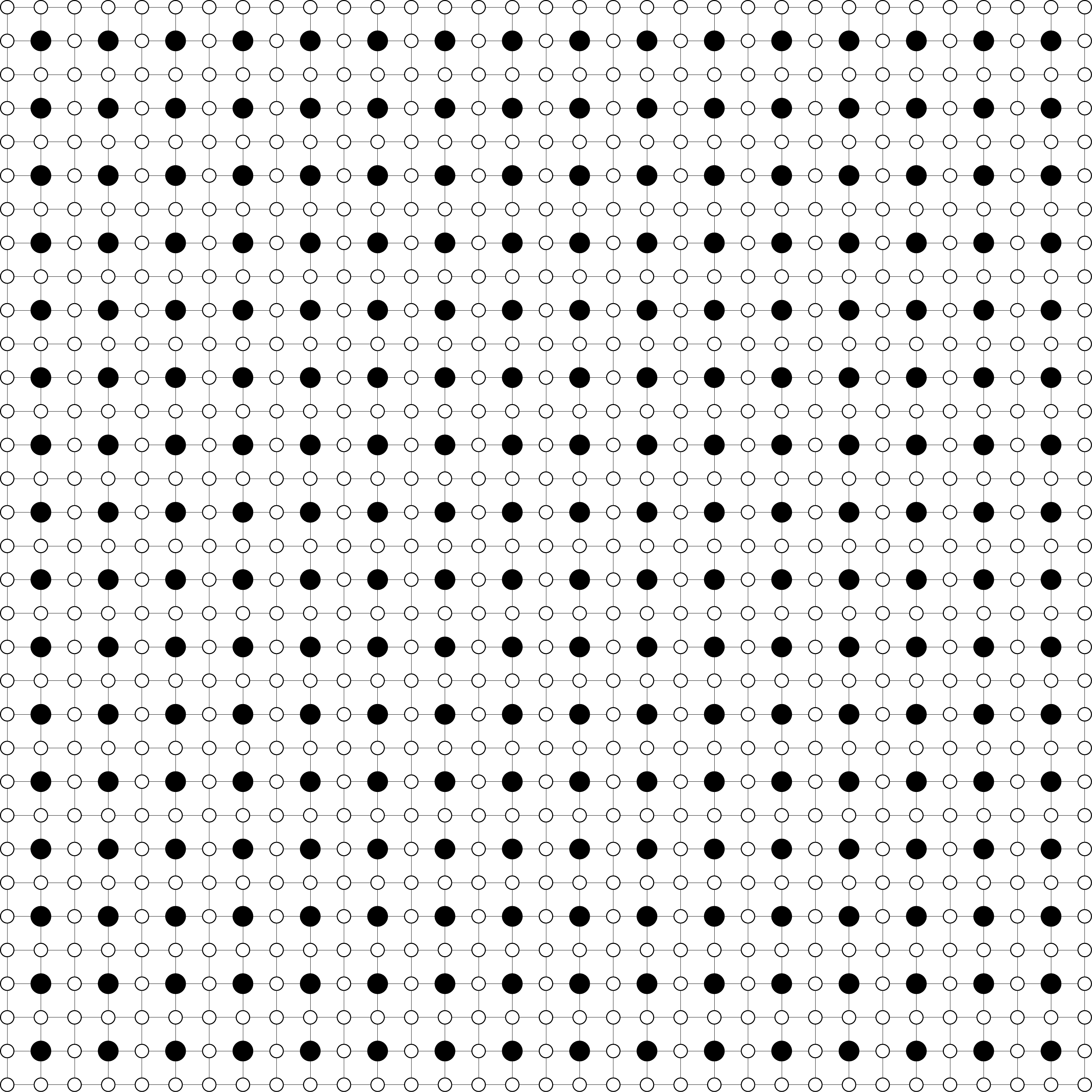}
                {\scriptsize \begin{align*} |\operatorname{det}(P_{c})| & = 4.228\cdot 10^{-76}\\ \kappa(P_{c}) &= \textbf{1.496}\end{align*}}%
            \end{minipage}
            \hfill 
            \begin{minipage}[b]{.36\textwidth}
                \hfill\includegraphics[width=1\textwidth]{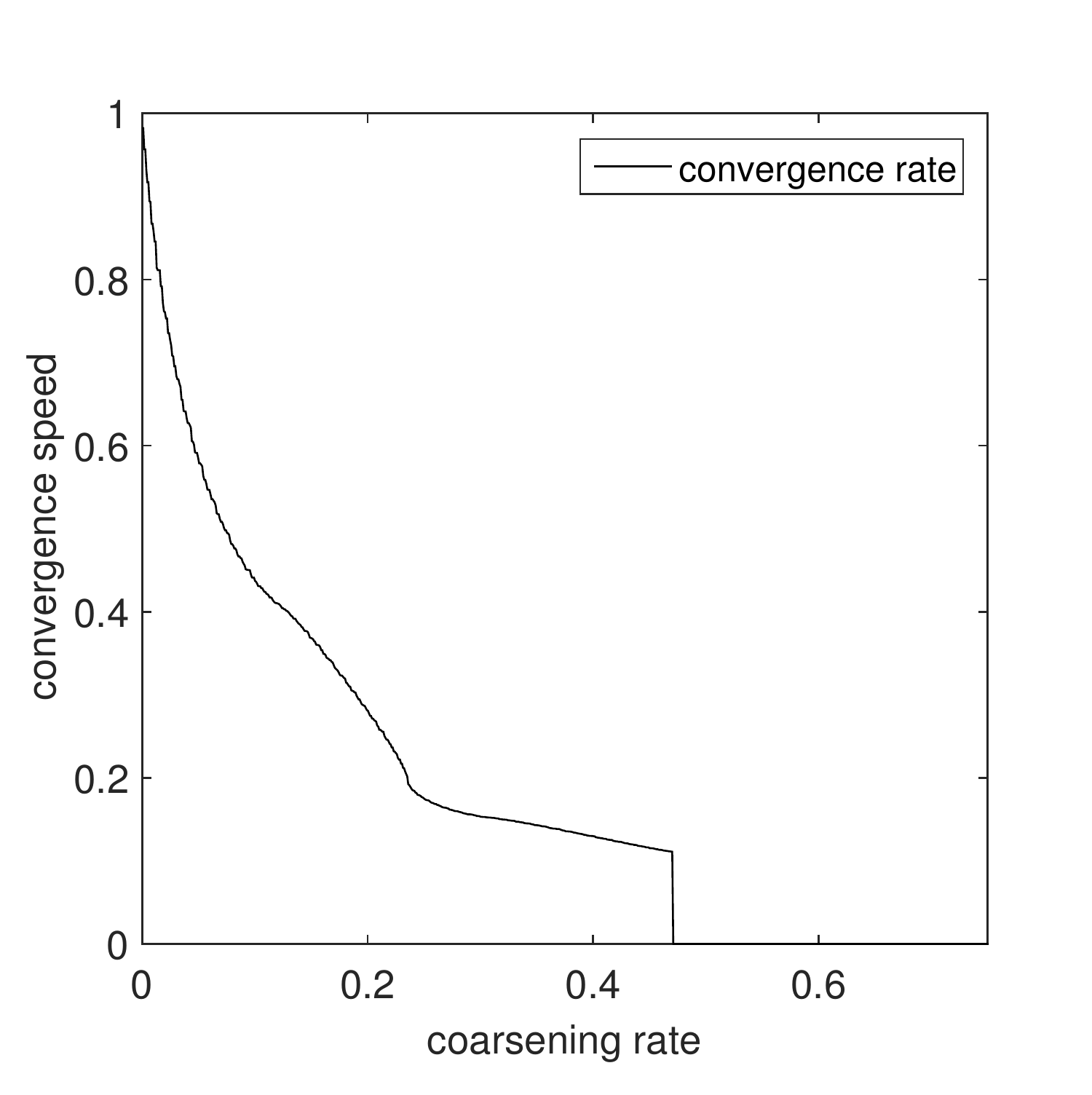}
            \end{minipage}
            \hfill \\
            \includegraphics[width=\textwidth]{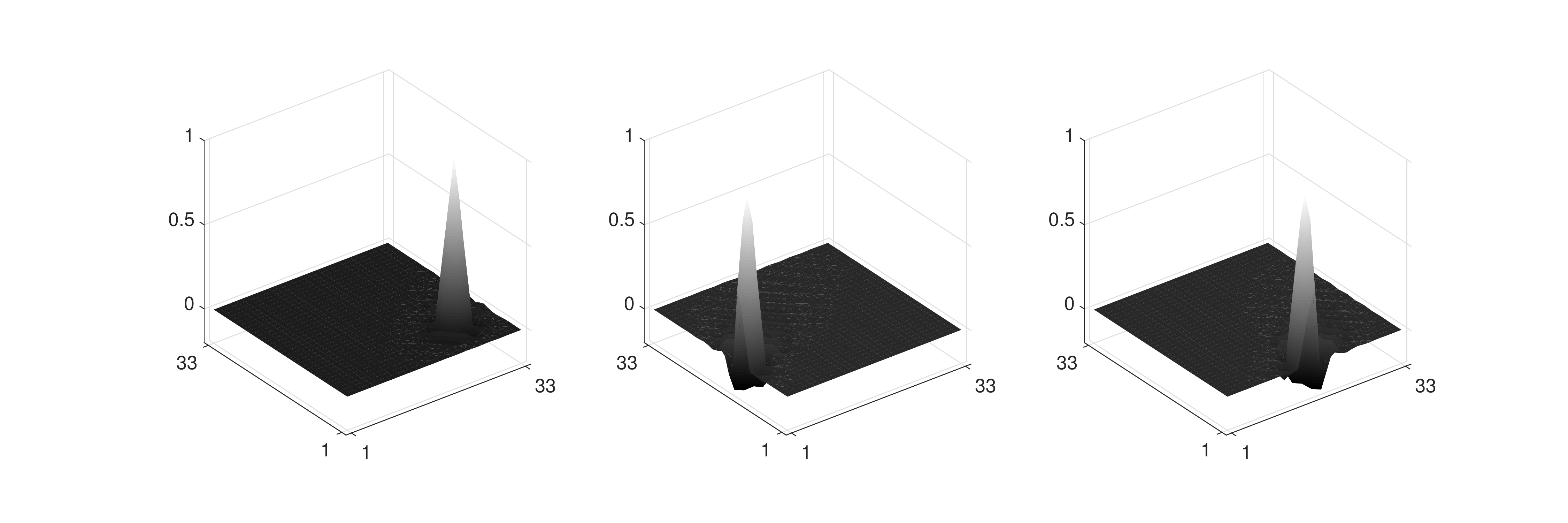} 
        \end{center}
        \caption{Lexicographic Gauss-Seidel for $P1$ ($k = 4$, $35\times 35$, $n_c = 289$). (top left) initial random choice of $C$; (top center) $C$ determined by \texttt{maxvol}; (top right) CR rate $\rho_{\Sopt} = 1 - \lambda_{n_c +1}$ w.r.t.\ ratio $\tfrac{n_{c}}{n}$; (bottom row) columns of $\PoptP$ with $C$ from \texttt{maxvol}}
        \label{fig:2d-spec}
    \end{figure}
    
    \begin{figure}[ht!]
        \begin{center}
            \hfill
            \begin{minipage}[b]{.25\textwidth}
                \includegraphics[width=1\textwidth]{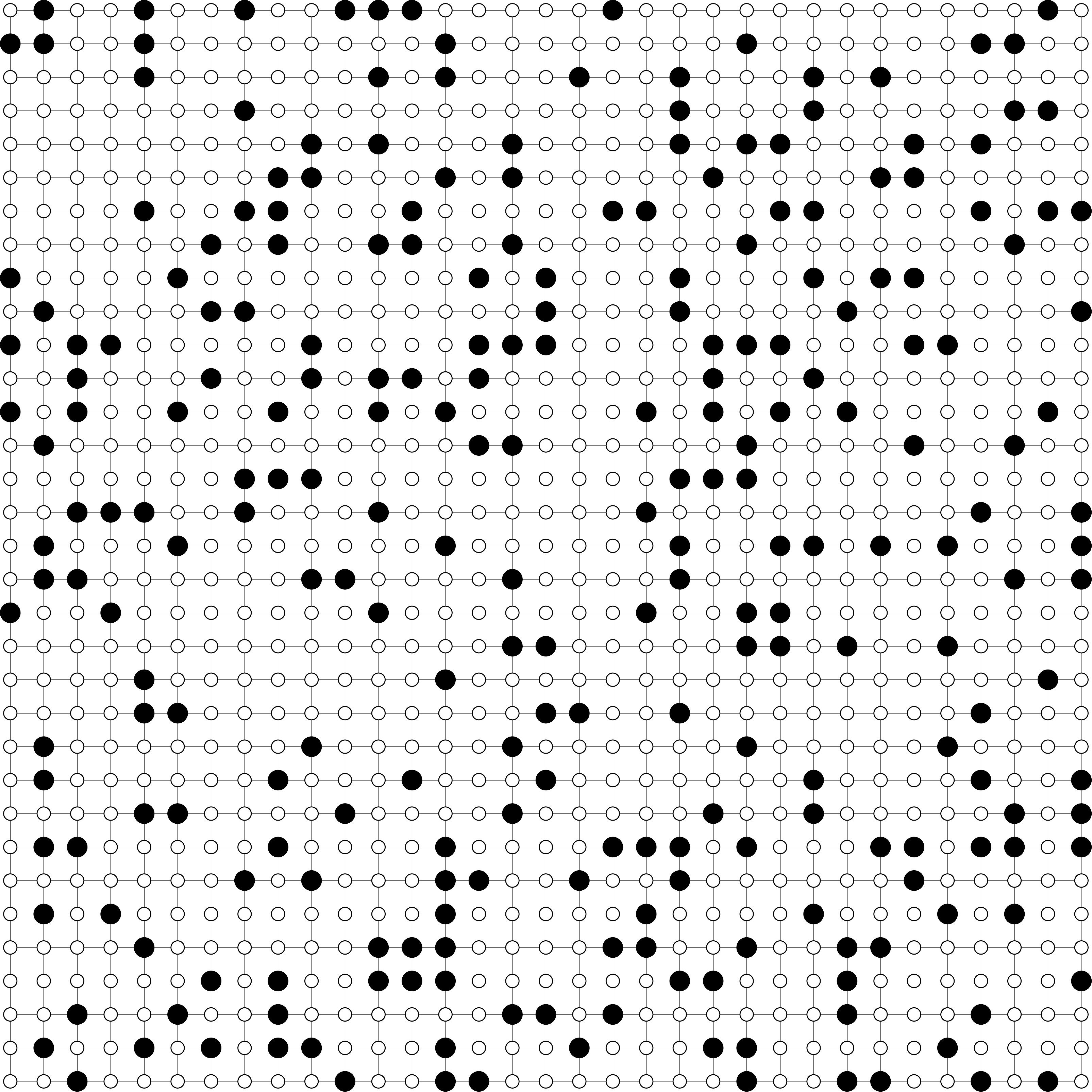}
                {\scriptsize \begin{align*} |\operatorname{det}(P_{c})| & = 2.071\cdot 10^{-240}\\ \kappa(P_{c}) &= 2.311\cdot 10^{11}\end{align*}}%
            \end{minipage}
            \hfill 
            \begin{minipage}[b]{.25\textwidth}
                \includegraphics[width=1\textwidth]{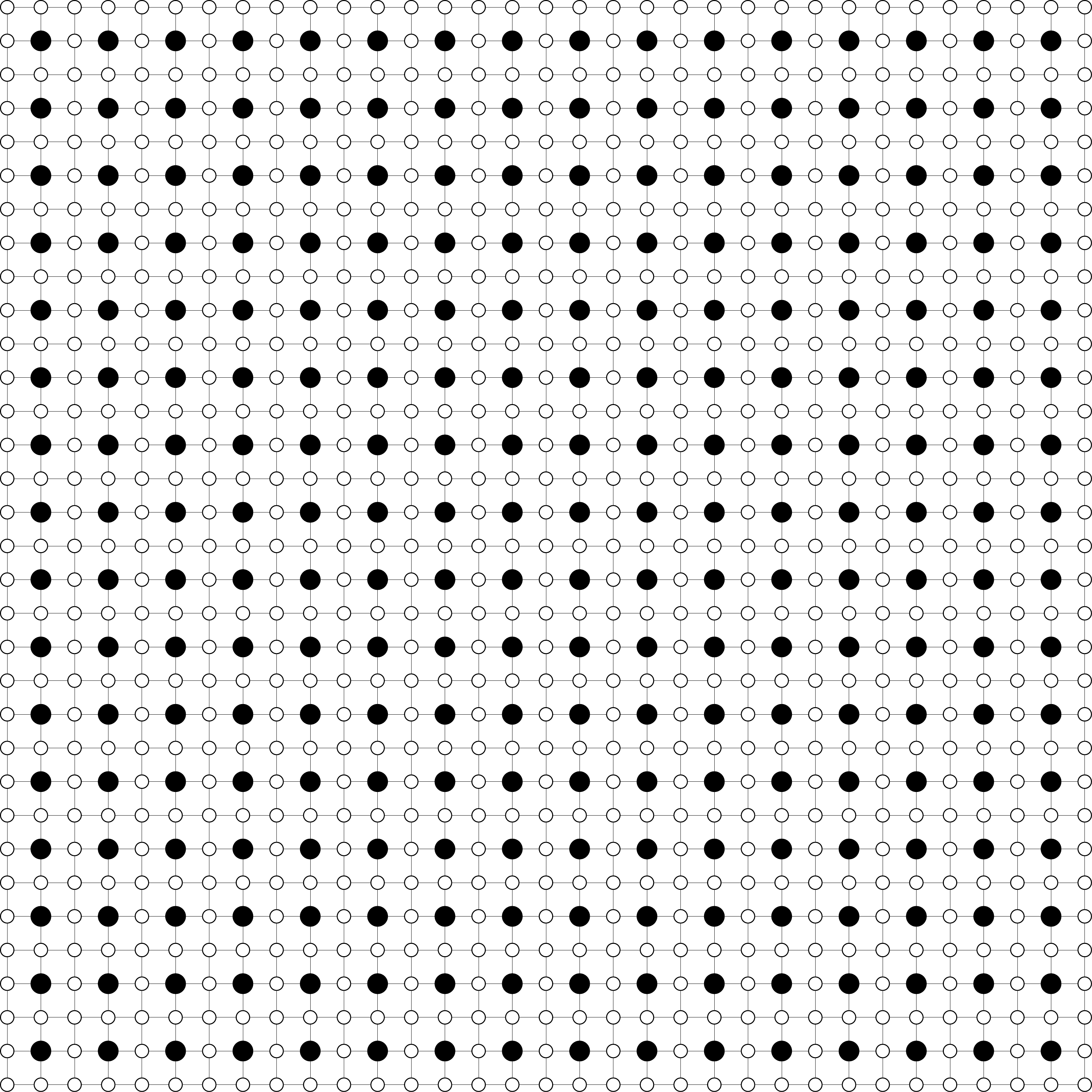}
                {\scriptsize \begin{align*} |\operatorname{det}(P_{c})| & = 1.545\cdot 10^{-76}\\ \kappa(P_{c}) &= \textbf{1.288}\end{align*}}%
            \end{minipage}
            \hfill 
            \begin{minipage}[b]{.36\textwidth}
                \hfill\includegraphics[width=1\textwidth]{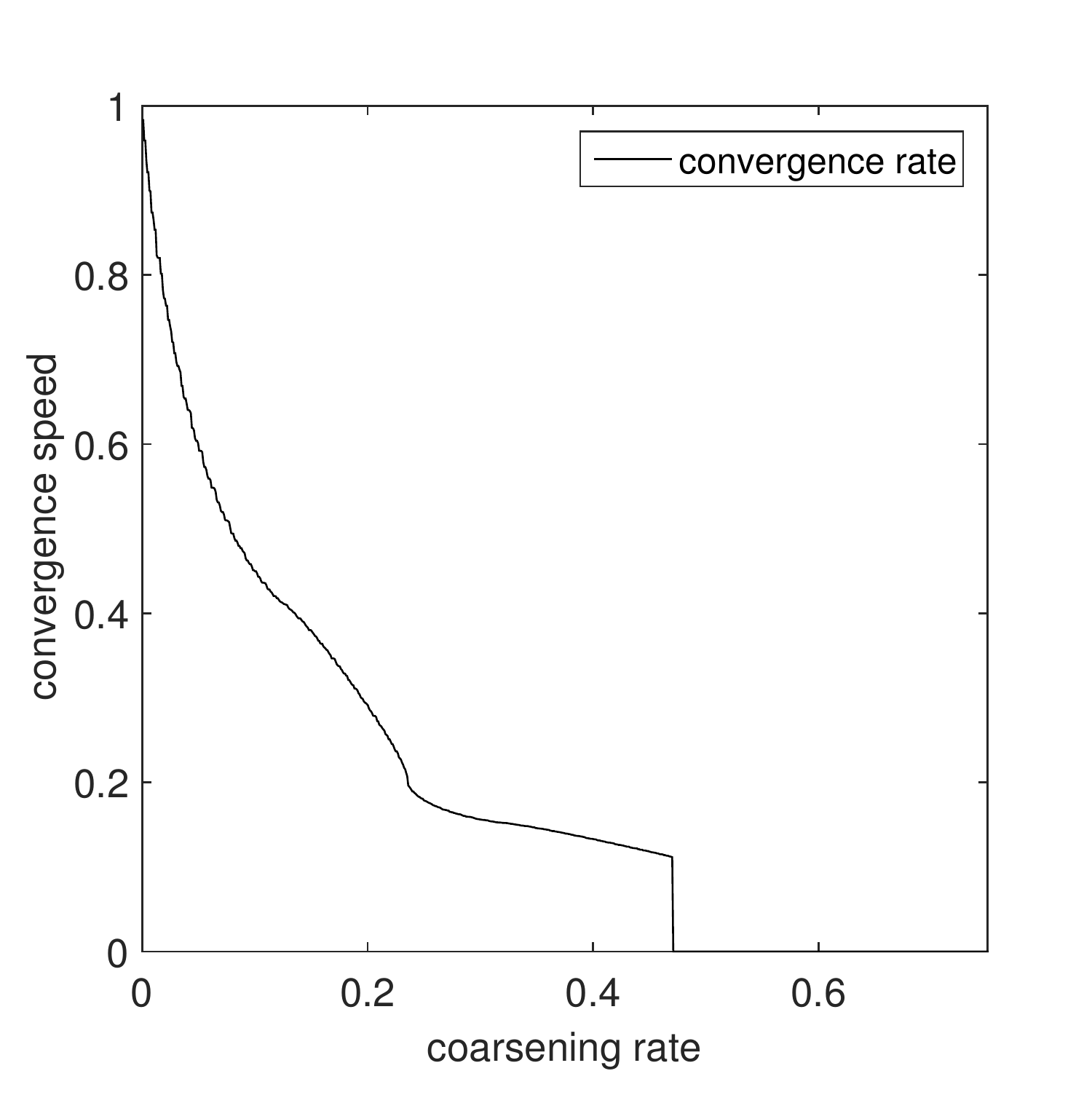}
            \end{minipage}
            \hfill \\
            \includegraphics[width=\textwidth]{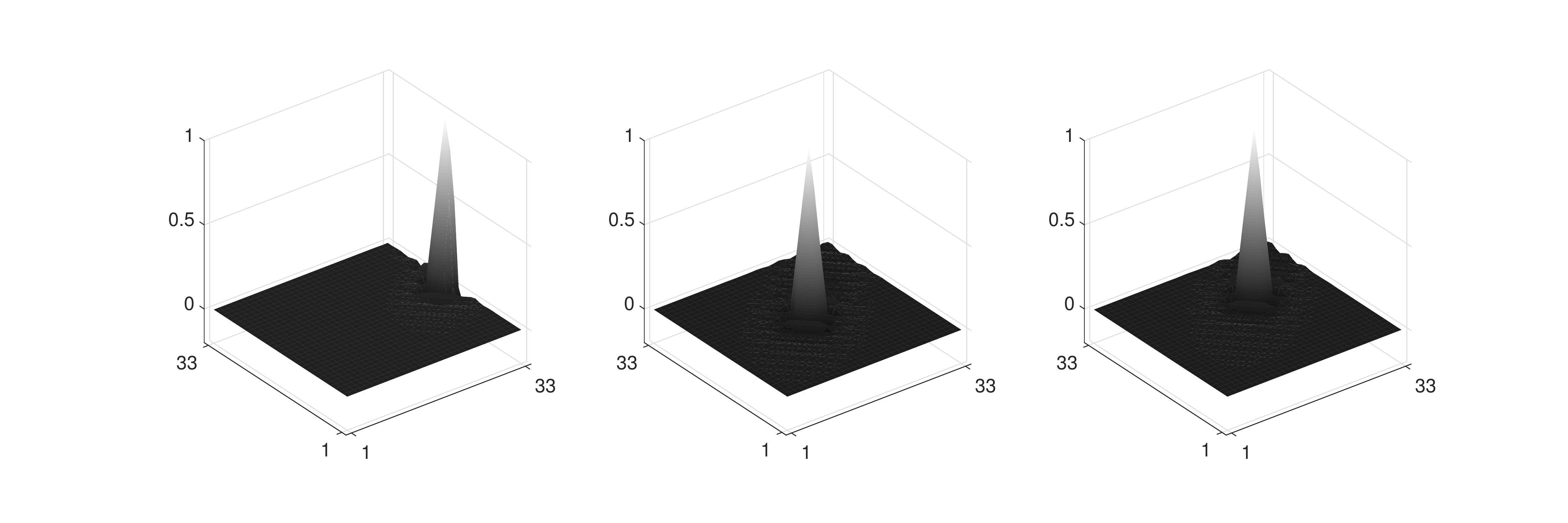} 
        \end{center}
        \caption{Lexicographic Gauss-Seidel for $P2$ ($k = 4$, $35\times 35$, $n_c = 289$). (top left) initial random choice of $C$; (top center) $C$ determined by \texttt{maxvol}; (top right) CR rate $\rho_{\Sopt} = 1 - \lambda_{n_c +1}$ w.r.t.\ ratio $\tfrac{n_{c}}{n}$; (bottom row) columns of $\PoptP$ with $C$ from \texttt{maxvol}}
        \label{fig:2d-spec1}
    \end{figure}
    \begin{figure}[ht!]
        \begin{center}
            \hfill
            \begin{minipage}[b]{.25\textwidth}
                \includegraphics[width=1\textwidth]{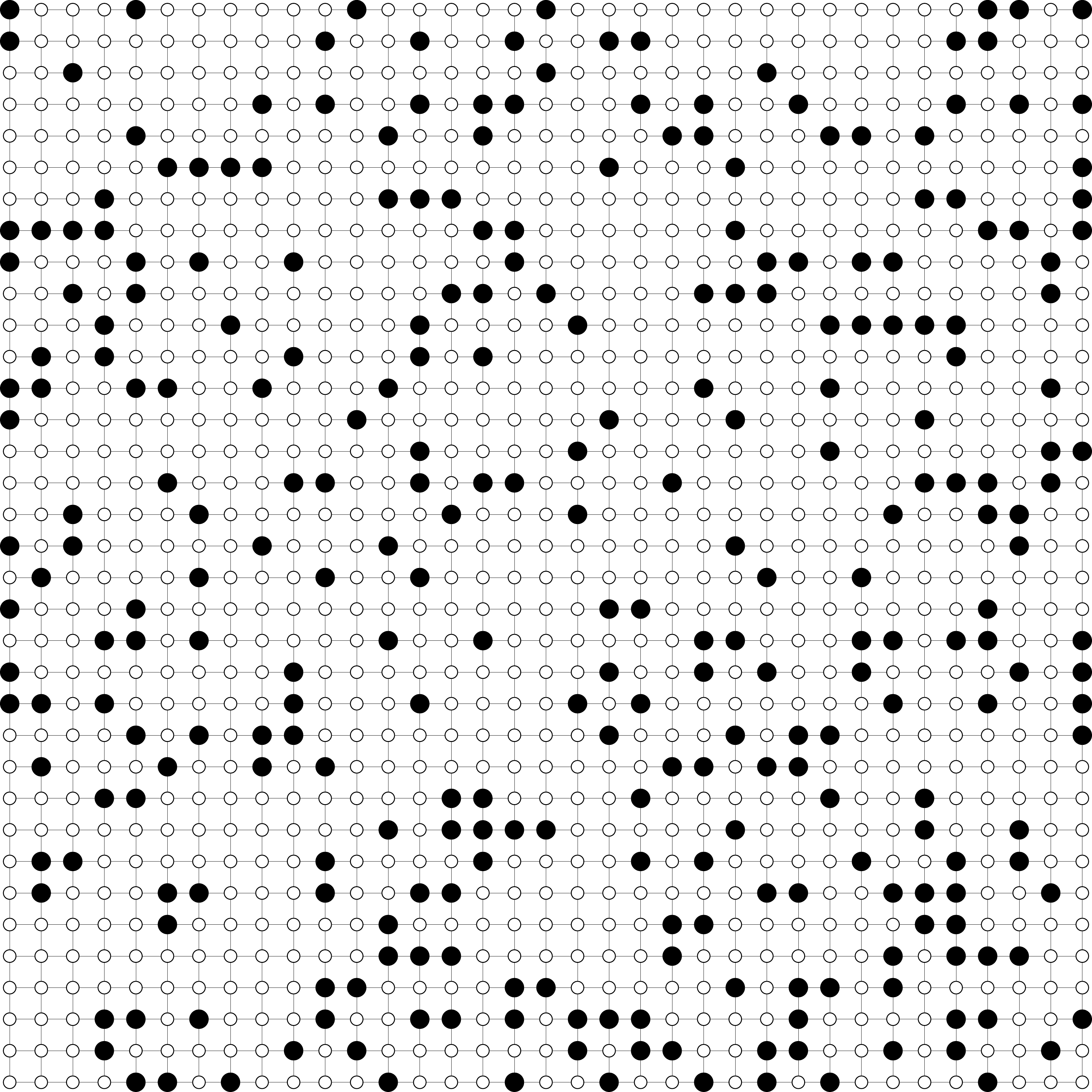}
                {\scriptsize \begin{align*} |\operatorname{det}(P_{c})| & = 4.888\cdot 10^{-255}\\ \kappa(P_{c}) &= 1.787128\cdot 10^{13}\end{align*}}%
            \end{minipage}
            \hfill 
            \begin{minipage}[b]{.25\textwidth}
                \includegraphics[width=1\textwidth]{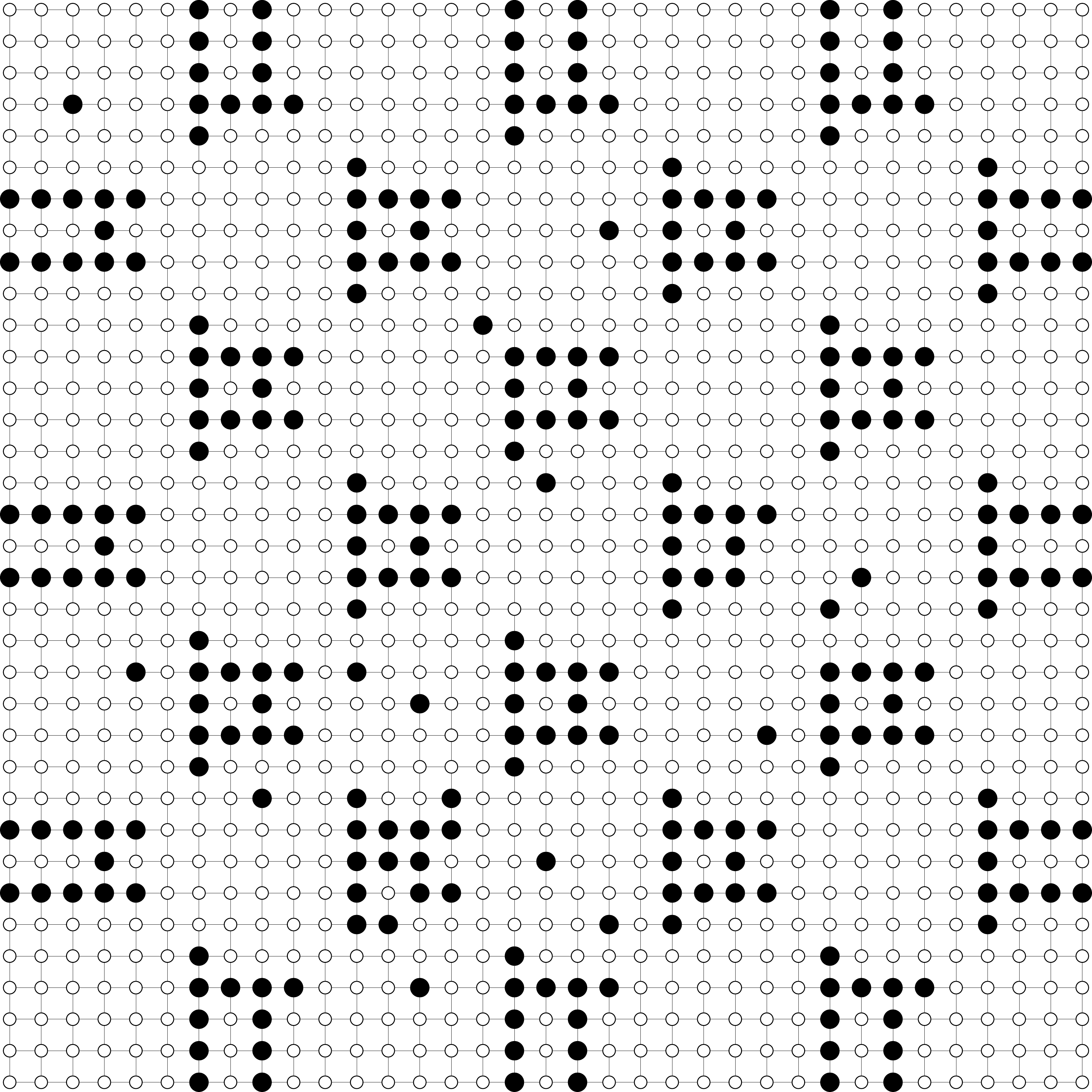}
                {\scriptsize \begin{align*} |\operatorname{det}(P_{c})| & = 8.594\cdot 10^{-47}\\ \kappa(P_{c}) &= \textbf{20.529}\end{align*}}%
            \end{minipage}
            \hfill 
            \begin{minipage}[b]{.36\textwidth}
                \hfill\includegraphics[width=1\textwidth]{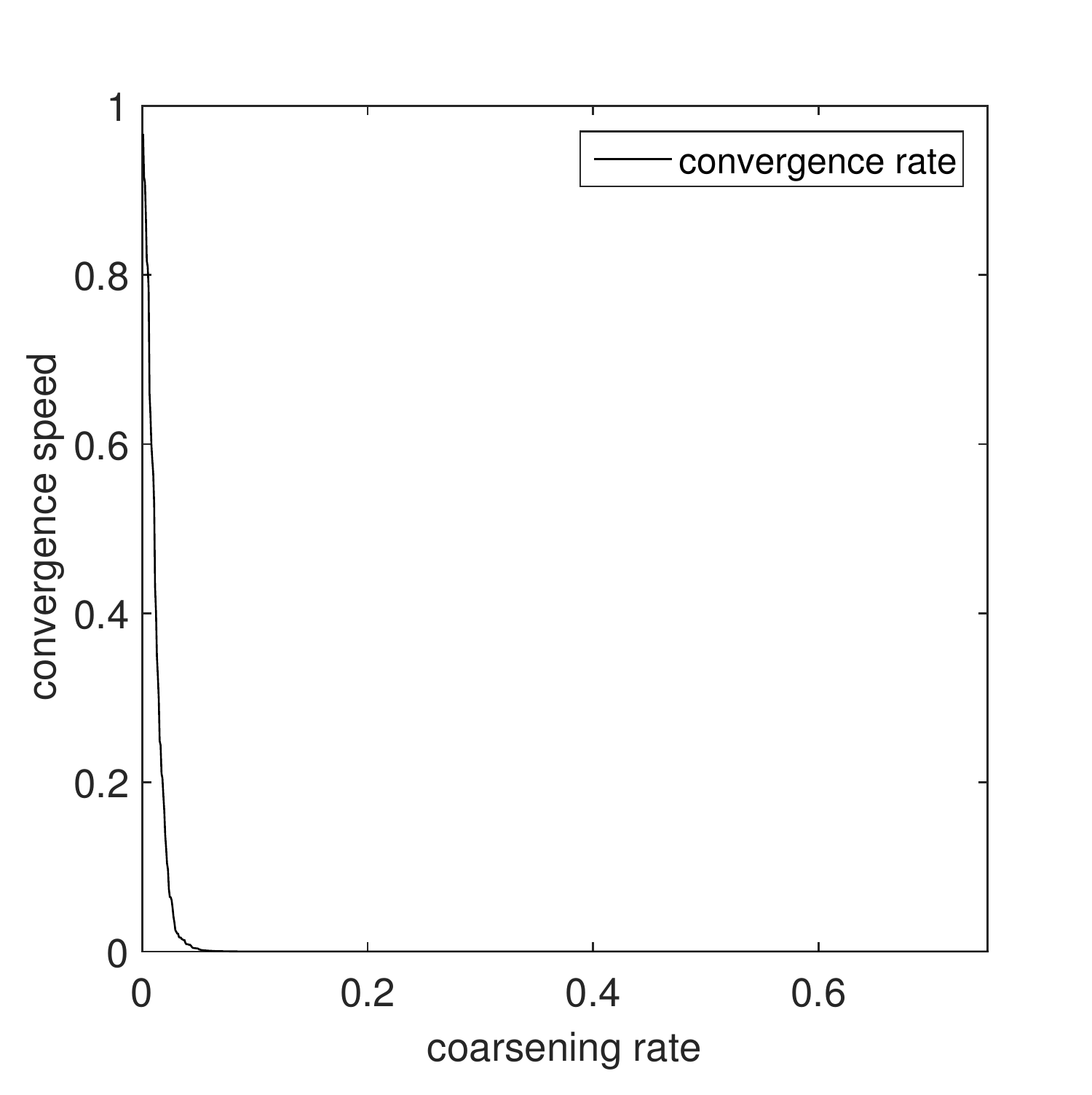}
            \end{minipage}
            \hfill \\
            \includegraphics[width=\textwidth]{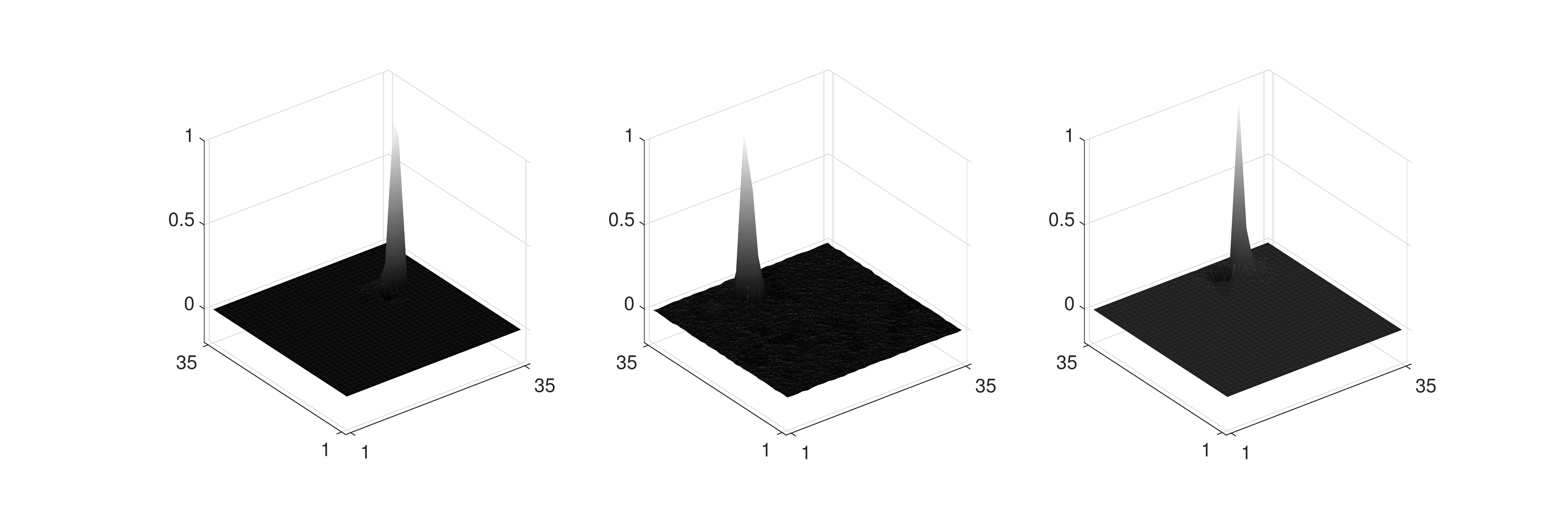} 
        \end{center}
        \caption{Red-black block Gauss-Seidel for $P4$, $5\times 5$ blocks ($k = 4$, $35\times 35$, $n_c = 144$). (top left) initial random choice of $C$; (top center) $C$ determined by \texttt{maxvol}; (top right) CR convergence rate $\rho_{\Sopt} = 1 - \lambda_{n_c +1}$  w.r.t.\ ratio $\tfrac{n_{c}}{n}$; (bottom row) three columns of $\PoptP$ with $C$ from \texttt{maxvol}}
        \label{fig:2d-spec3}
     \end{figure}

\subsection{Bootstrap AMG and the generalized eigenproblem}
In this section, we consider designing a bootstrap AMG setup algorithm
that aims at solving the generalized eigen-problem in~\eqref{eq:gep} to compute test vectors used in 
constructing least squares interpolation. We note that
in the original BAMG setup algorithm developed in~\cite{BAMG2010} the eigen-problem involving 
only $A$ was used to compute approximate eigenvectors that are provided as input
to the least squares process that constructs the interpolation matrix. As we will point out there are
only a few changes that are needed in order to adapt the least squares interpolation and bootstrap multilevel setup
to the generalized eigenvalue problem.

Plots of the eigenvalues of $A$ and $(A,\XMt)$ for various tests with $N=16$ are provided in Figure~\ref{fig:spec}.
These results show that for Problem P4 the spectrum of $A$ and 
$(A,\XMt)$ differ substantially, especially with respect to the number of near-null eigenvectors.
This observation together with the fact that the optimal
interpolation matrix has as columns eigenvectors of the generalized eigen-problem involving
$A$ and $\XMt$ suggest that AMG interpolation should be based on approximating these
generalized eigenvectors. To test this idea we compare the BAMG setup that uses eigen-approximations
of $A$ to the one that uses eigen-approximations of $\XMt^{-1}A$.

\begin{figure}[ht!]
    \hfill
    \begin{minipage}[b]{.3\textwidth}
        \includegraphics[width=1\textwidth]{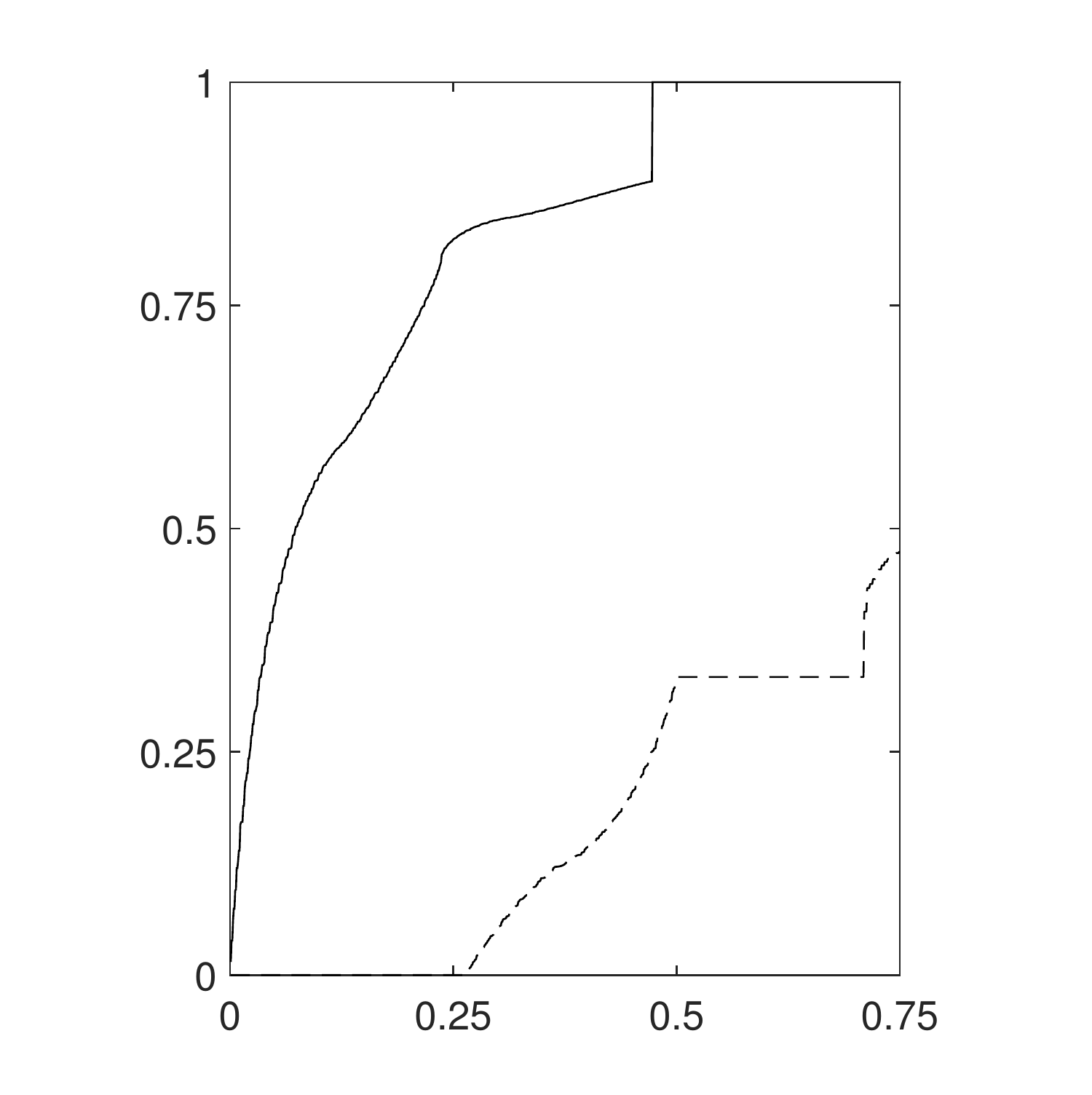}
    \end{minipage}
    \hfill 
    \begin{minipage}[b]{.3\textwidth}
        \includegraphics[width=1\textwidth]{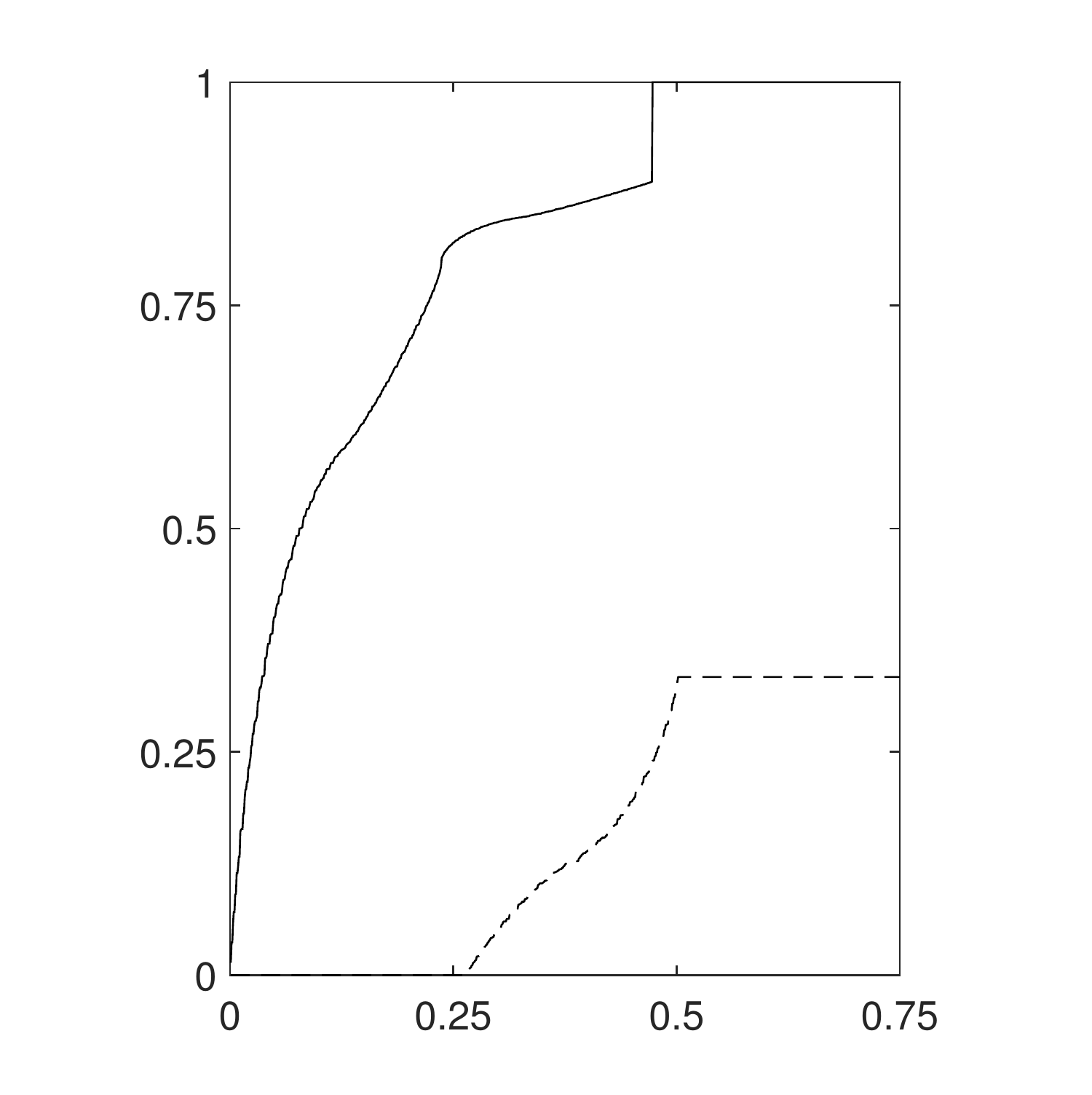}       
    \end{minipage}
    \hfill 
    \begin{minipage}[b]{.3\textwidth}
        \includegraphics[width=1\textwidth]{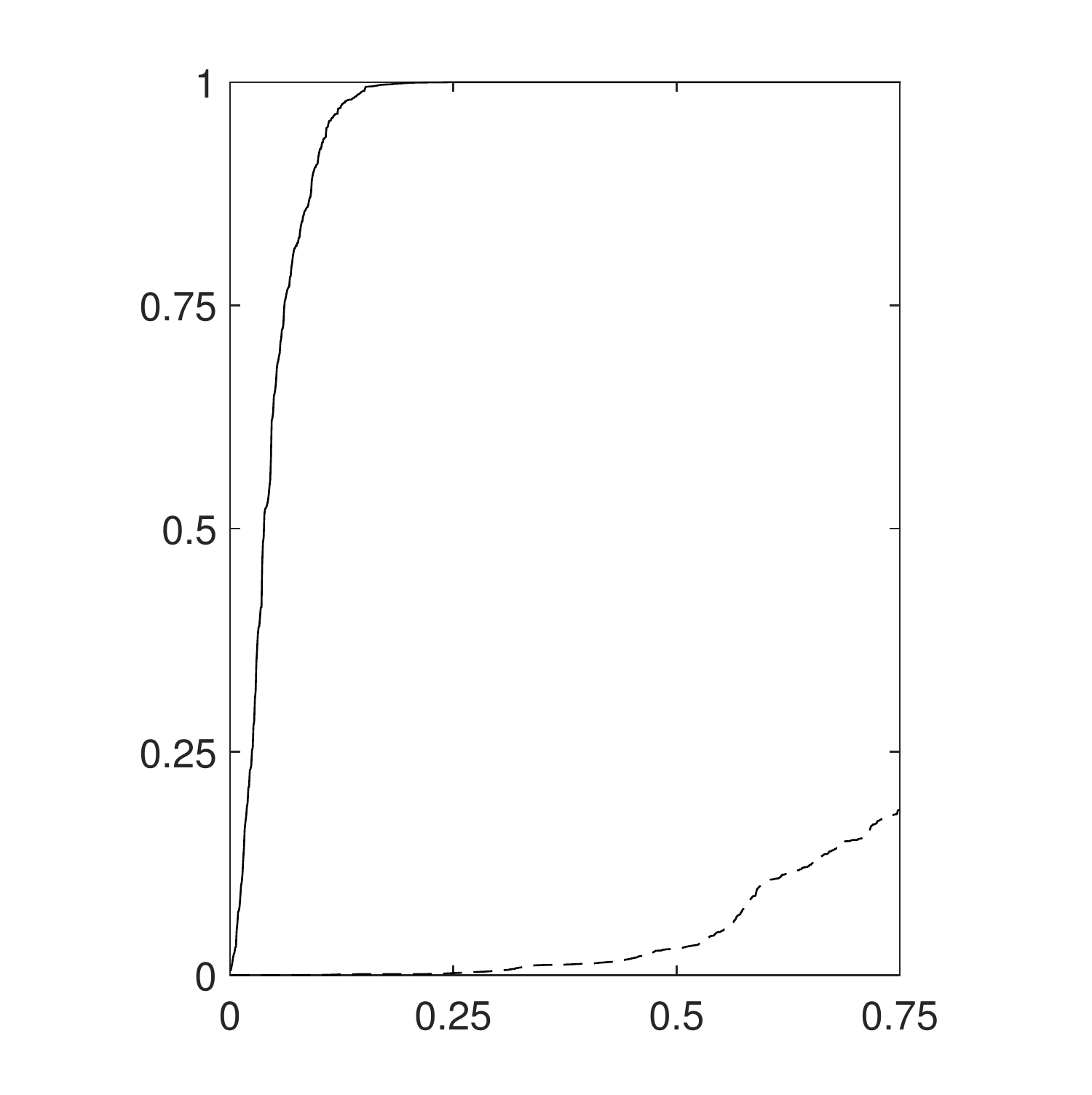}
    \end{minipage}
    \hfill
    \caption{Lower $\tfrac{3}{4}$ of the spectra of $(A,\XMt)$ (solid line) and $A$ (dashed line), which is scaled by $\lambda_{\rm max}(A)$. (left) Lex.\ Gauss-Seidel for P1 ($k=4$, $35 \times 35$); (center) Lex.\ Gauss-Seidel for P2 ($k=4$, $35 \times 35$); (right) Red-black block Gauss-Seidel with $5\times 5$ blocks for P4 ($k=4$, $35 \times 35$).
}    \label{fig:spec}
\end{figure}
\subsubsection{Bootstrap AMG setup}
For definiteness, we provide an outline of the version of the bootstrap AMG setup that was 
designed and analyzed in~\cite{BAMG2010} and that we use in our tests. We point out the modifications that are necessary to adjust
the components of the method to the generalized eigenvalue problem.

The least squares interpolation matrix used in BAMG is defined to fit
collectively a set of test vectors (TVs) that should characterize the eigenvectors with small eigenvectors of the 
system matrix (or the generalized eigenproblem).   
Assuming the sets of interpolatory variables, $C_i$, for each $i\in F$, and a set of test vectors, $\tvV = \{v^{(1)}, \ldots,
v^{(k)}\}$, have been determined, the $i$th row of $P$, denoted by $p_i$,
is defined as the minimizer of the local least squares problem: 
\begin{equation}\label{eq:LSfuncrowi}
\mathcal{L}(p_i) = 
\sum_{\kappa=1}^k\omega_k\left(v_{\{i\}}^{(\kappa)} - \sum_{j\in C_{i}} \left(p_{i}\right)_{j} v_{\{j\}}^{(\kappa)}\right)^{2} \rightarrow \min.
\end{equation} 
Here, the notation $v_{\widetilde{\Omega}}$ denotes the canonical restriction of the vector $\widetilde{v}$ to
the set $\widetilde{\Omega} \subset \Omega$, e.g.,
$v_{\{i\}}$ is simply the $i$th entry of $v$. Conditions on the uniqueness of the solution to minimization problem
\eqref{eq:LSfuncrowi} and an explicit form of the minimizer have been
 derived in \cite{BAMG2010}. 
 In contrast to the original least squares interpolation we choose the weights $\omega_{\kappa}>0$ by
\[
\omega_{\kappa} = \frac{\norm[X]{v^{(\kappa)}}}{\norm[A]{v^{(\kappa)}}}.
\] In the original formulation $X=I$ corresponds to the inverse Rayleigh-Quotient of $v^{(\kappa)}$. Choosing $X = \XMt$ the weight is 
the inverse generalized Rayleigh-Quotient w.r.t.~the pencil $(A,\XMt)$.

\begin{figure}
    \begin{tikzpicture}[scale=1.5]
        \draw[thin] (-.25,-.25) grid (4.25,2.25);
        \foreach \x in {0,1,...,4}{
            \foreach \y in {0,1,2}{
                \node[draw=black,fill=black!50!white,circle,inner sep=0cm,minimum width=.2cm] (A\x\y) at (\x,\y) {} ;
            }
        }
        
        \foreach \x in {0,2,...,4}{
            \foreach \y in {0,2}{
                \node[draw=black,fill=white,circle,inner sep=0cm,minimum width=.2cm] (C\x\y) at (\x,\y) {} ;
            }
        }
        
        \path[-latex,black,draw,shorten >=2pt] (A00) to[bend right] (A10);
        \path[-latex,black,draw,shorten >=2pt] (A20) to[bend right] (A10);
        
        \path[-latex,black,draw,shorten >=2pt] (A00) to[bend right] (A01);
        \path[-latex,black,draw,shorten >=2pt] (A02) to[bend right] (A01);
        
        \path[-latex,black,draw,shorten >=2pt] (A42) to (A31);
        \path[-latex,black,draw,shorten >=2pt] (A22) to (A31);
        \path[-latex,black,draw,shorten >=2pt] (A40) to (A31);
        \path[-latex,black,draw,shorten >=2pt] (A20) to (A31);
    \end{tikzpicture}
    \caption[Full coarsening and interpolation relations]{Full coarsening and interpolation relations (\raisebox{1pt}{\resizebox{!}{.5em}{\begin{tikzpicture}\node[draw=black,fill=black!50!white,circle,inner sep=0cm,minimum width=.2cm]{};\end{tikzpicture}}} $F$, \raisebox{1pt}{\resizebox{!}{.5em}{\begin{tikzpicture}\node[draw=black,fill=white,circle,inner sep=0cm,minimum width=.2cm]{};\end{tikzpicture}}} $C$).} \label{fig:fc}
\end{figure}

The test vectors used in the least squares process are computed using a bootstrap multilevel setup cycle. The algorithm begins with 
relaxation applied to the homogenous system, 
\begin{equation}\label{homog}
A_lx_l =0,
\end{equation}
on each grid, $l = 0,..., L-1$; assuming that a priori knowledge of the 
algebraically smooth error is not available, these vectors are initialized randomly on the finest grid, whereas on all coarser grids
they are defined by restricting the existing test vectors computed on the previous finer grid.  

Once an initial MG hierarchy has been computed, 
the current sets of TVs are further enhanced on all grids using the 
existing multigrid structure.  Specifically, the given hierarchy is used to formulate a multigrid eigensolver
which is then applied to an appropriately chosen generalized eigenproblem to compute additional test vectors.  
This overall process is then repeated with the current AMG method applied in addition to (or replacing)
relaxation as the solver for the homogenous systems in \eqref{homog}.
Figure \ref{fig:boot:setupcycle} provides an schematic outline of the bootstrap $V$- and $W$-cycle setup algorithms.  {In general, $V^m$ and $W^m$ denote setup algorithms that use $m$ iterations of the $V$- and $W$-cycles, respectively.}

The rationale behind the multilevel generalized eigensolver (MGE) is as follows.  Assume an initial multigrid hierarchy has been constructed. 
Given the initial Galerkin operators $A_{0}, A_{1}, \ldots, A_{L}$ on each level
and the corresponding interpolation operators $P_{l+1}^{l}, l =
0,\ldots,L-1$, define the composite interpolation operators as
$P_{l} = P_{1}^{0}\cdot \ldots \cdot P_{l}^{l-1},\ l = 1, \ldots, L$.  Then, 
for any given vector $x_{l} \in \mathbb{R}^{n_l}$ we have
$
  \innerprod[A_{l}]{x_{l}}{x_{l}} = \innerprod[A]{P_{l}x_{l}}{P_{l}x_{l}}.
$  
Furthermore, defining $X_{l} = P_{l}^{T}X P_{l}$ for any $X$ symmetric and positive definite we obtain
\begin{equation*}
  \frac{\innerprod[A_{l}]{x_{l}}{x_{l}}}{\innerprod[X_{l}]{x_{l}}{x_{l}}} = \frac{\innerprod[A]{P_{l}x_{l}}{P_{l}x_{l}}}{\innerprod[X]{P_{l}x_{l}}{P_{l}x_{l}}}.
\end{equation*} 
This observation in turn implies that 
on any level $l$, given a vector $x^{(l)} \in \mathbb{R}^{n_l}$ and $\lambda^{(l)}\in
  \mathbb{R}$ such that
$A_l x^{(l)} = \lambda^{(l)}X_l x^{(l)},$
its Rayleigh quotient (RQ) fullfills
\begin{equation}\label{eq:evalapprox}
  \mbox{rq}(P_l x^{(l)}) :=
  \frac{\innerprod[A]{P_{l}x^{(l)}}{P_{l}x^{(l)}}}{\innerprod[X]{P_{l}x^{(l)}}{P_{l}x^{(l)}}}
  = \lambda^{(l)}.
\end{equation} This provides a relation among the eigenvectors and eigenvalues
of the operators in the multigrid hierarchy on all levels with 
the eigenvectors and eigenvalues of the finest-grid matrix pencil $(A,X)$. 
Again one obtains the original formulation of the bootstrap setup cycle by choosing 
$X = I$ and by choosing $X=\XMt$ one obtains a setup cycle that yields approximations to the eigenvectors of the
matrix pencil $(A,\XMt)$.
We note that the eigenvalue approximations in~\eqref{eq:evalapprox} are continuously updated within the algorithm so that the overall approach resembles an inverse Rayleigh-Quotient
iteration found in eigenvalue computations (cf.~\cite{JWilkinson_1965}).  
For additional details of the algorithm and its implementation we refer to the paper~\cite{BAMG2010}.  

\begin{figure}
\begin{center}
     \tikzstyle{greenpoint}=[circle,inner sep=0pt,minimum size=2mm,draw=black!100,fill=black!30]
    \tikzstyle{whitepoint}=[rectangle,inner sep=0pt,minimum size=2mm,draw=black!100,fill=black!80]
    \tikzstyle{redpoint}=[diamond,inner sep=0pt,minimum size=2.55mm,draw=black!100,fill=black!80]
    \tikzstyle{blackpoint}=[circle,inner sep=0pt,minimum size=2mm,draw=black!100,fill=black!100]
    \tikzstyle{bluepoint}=[circle,inner sep=0pt,minimum size=2mm,draw=black!100,fill=black!0]
    \resizebox{\textwidth}{!}{\begin{tikzpicture}
      \draw [sharp corners] (0,0) node[blackpoint] {} --
      ++(300:1cm) node[blackpoint] {} --
      ++(300:1cm) node[blackpoint] {} --
      ++(300:1cm) node[blackpoint] {} --
      ++(300:1cm) node[redpoint] (L4) {} --
      ++(60:1cm) node[bluepoint] (L3) {} --
      ++(60:1cm) node[bluepoint] (L2) {} --
      ++(60:1cm) node[bluepoint] (L1) {} --
      ++(60:1cm) node[whitepoint] (L0) {} --
      ++(0:1cm) node[greenpoint] {} --
      ++(300:1cm) node[greenpoint] (end) {};
      \draw[dashed] (end) --  ++(300:1cm);
      
      \draw (L0) + (3cm,0) node[blackpoint,label=right:{\small Relax on $Av=0, v \in \mathcal{V}^r$, compute $P$}] {};
      \draw (L1) + (3.5cm,0) node[redpoint,label=right:{\small Compute $v$, s.t., $Av=\lambda Tv$, update $\mathcal{V}^e$}] {};
      \draw (L2) + (4.0cm,0) node[bluepoint,label=right:{\small Relax on $Av = \lambda T v, v \in \mathcal{V}^e$}] {};
      \draw (L3) + (4.5cm,0) node[greenpoint,label=right:{\small Relax on and / or solve using MG $Av=0, v \in \mathcal{V}^r$ and $Av = \lambda T v, v \in \mathcal{V}^e$, recompute $P$}] {};
      \draw (L4) + (5cm,0) node[whitepoint,label=right:{\small Test
        MG method, update $\mathcal{V}$}] {};

      \draw (L3) +(0,4cm) node (start) {};
      \draw[sharp corners] (start) ++(300:-3cm) node[blackpoint] {} --
      ++(300:1cm) node[blackpoint] {} --
      ++(300:1cm) node[blackpoint] {} --
      ++(300:1cm) node[redpoint] {} --
      ++(60:1cm) node[greenpoint] {} --
      ++(300:1cm) node[redpoint] {} --
      ++(60:1cm) node[bluepoint] {} --
      ++(60:1cm) node[greenpoint] {} --
      ++(300:1cm) node[greenpoint] {} --
      ++(300:1cm) node[redpoint] {} --
      ++(60:1cm) node[greenpoint] {} --
      ++(300:1cm) node[redpoint] {} --
      ++(60:1cm) node[bluepoint] {} --
      ++(60:1cm) node[greenpoint] {} --
      ++(300:1cm) node[greenpoint] {} --
      ++(300:1cm) node[redpoint] {} --
      ++(60:1cm) node[bluepoint] {} --
      ++(60:1cm) node[bluepoint] {} --
      ++(60:1cm) node[whitepoint] {} -- 
      ++(0:1cm) node[greenpoint] {} --     
      ++(300:1cm) node[greenpoint] {} --
      ++(300:1cm) node[greenpoint] {} --
      ++(300:1cm) node[redpoint] {} --
      ++(60:1cm) node[greenpoint] (end) {};
      \draw[dashed] (end) --  ++(300:1cm);
    \end{tikzpicture}}
  \caption{Galerkin Bootstrap AMG $W$-cycle and $V$-cycle setup schemes.
  Here, $\mathcal{V}^r$ denotes the set of relaxed vectors and $\mathcal{V}^e$ the set of vectors coming from the MGE process.
  \label{fig:boot:setupcycle}}
\end{center}
\end{figure}
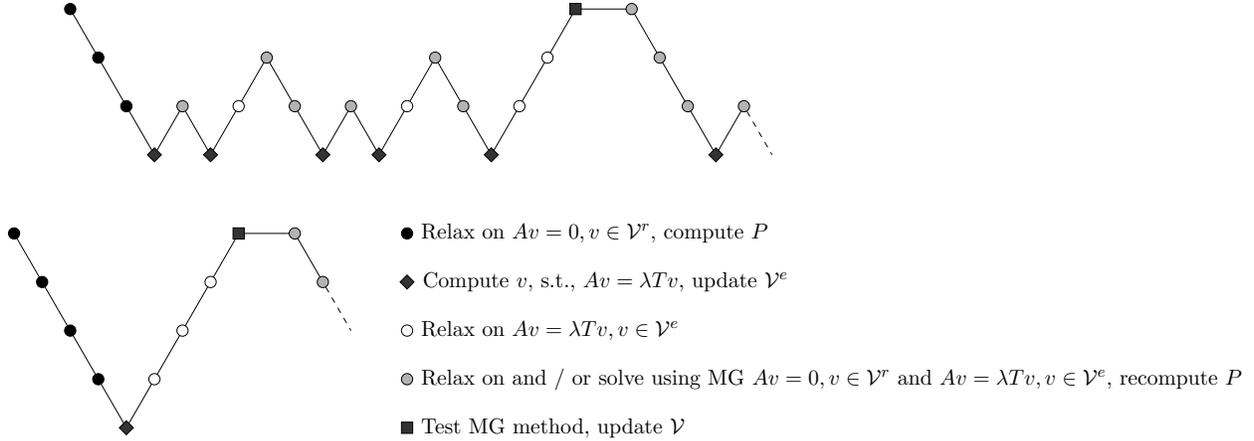

To illustrate the effect of the choice of $T_l$ in the MGE bootstrap process, we provide results in which a $W^2$-cycle bootstrap setup using four forward Gauss Seidel pre- and four backward Gauss Seidel post-smoothing steps to compute the set of relaxed vectors $\mathcal{V}^r$ and set of bootstrap vectors $\mathcal{V}^e$ coming from the MGE process, with $k_r=|\mathcal{V}^r| = 8$ and $k_e= |\mathcal{V}^e|=8$, respectively.  The sets $\mathcal{V}^r$ and $\mathcal{V}^e$ are then combined to form the set of TVs $\mathcal{V}$ that is used to compute the least squares interpolation operator on each level.   
In these tests, only relaxation is applied to the homogenous systems in both setup cycles to update the sets $ \mathcal{V}^r$.
We use forward Gauss Seidel as a pre-smoother and backward Gauss Seidel as the post-smoother in the solve phase as well.  
The coarse grids and sparsity structure of interpolation are defined as in the previous tests on all levels, i.e., by full coarsening and nearest
neighbor interpolation as depicted in Figure 4.3, and the problem is coarsened to a coarsest level with $h=1/8$.

The results of our experiments are given in Table~\ref{tab:bamg}.
Here, we observe that using the generalized eigen-problem
in the BAMG setup gives uniformly better results, than just working with eigen-approximations of $A$.  
In addition, if we compare the results obtained in the left table with the results given earlier
(see Table~\ref{tab:CR}), then we see that the multilevel BAMG approach (with sparse $P$) converges 
faster than the two-level method that uses the ideal interpolation operator. 
\begin{center}
	\begin{table}[ht!]
		\begin{tabular}{@{\extracolsep{-1.mm}} |c|cccc|}
			\hline   
			Size / k   & 1 & 2 & 4 & 8  \\ \hline
			$17^2$  & .276&.377 &.398&.626      \\             
			$33^2$  & .260&.256 &.302&.445     \\                  
			$65^2$  & .261&.256 &.299&.427       \\             
			$129^2$& .261&.256 &.299&.427      \\            
			\hline 
		\end{tabular}
	\hspace{.25cm}
		\begin{tabular}{@{\extracolsep{-1.mm}} |c|cccc|}
			\hline   
			Size / k   & 1 & 2 & 4 & 8  \\ \hline
			$17^2$  & .357& .592 &.405&.966      \\             
			$33^2$  & .416&.591 &.302&.953     \\                  
			$65^2$  & .261&.256 &.303&.573       \\             
			$129^2$& .261&.256 &.305&.571      \\            
			\hline 
		\end{tabular}
		\vspace{.3cm}
		\caption{BAMG setup results for Problem P4 with $X = \XMt$ (left) and $X = I$ (right).}\label{tab:bamg}
\end{table}
\end{center}
\vspace{-1.5cm}

\section{Conclusion}
In this paper, we introduced an optimal form of classical AMG interpolation and a measure of the quality of the coarse variable set that gives precise estimates of the convergence rate of two-grid method with this optimal choice of interpolation, which is based on eigenvectors of the 
generalized eigen-problem involving the system matrix and its associated 
symmetrized smoother.  We derived the equivalence between the ideal and optimal forms of interpolation in the case of reduction-based AMG
(i.e., when $F$-relaxation is used) and we showed numerically that in the case of full smoothing the convergence rates of 
two-grid methods with these different choices can vary substantially.  We also showed that for full smoothing and a proper choice of the coarse
variable type, the optimal and ideal interpolation matrices have the same range.  Finally, using these new results 
we designed a generalized bootstrap AMG setup algorithm that incorporates this generalized eigen-problem and illustrated the utility of the approach
when applied to a scalar diffusion problem with highly varying diffusion coefficient.  Numerically, we observed that the 
BAMG method (spanning multiple levels) with sparse $P$ outperforms the two-grid method with the ideal $P$ (which is
a dense matrix).  
In addition, in our tests of the BAMG algorithm we did not use strength of connection in defining the sparsity structure of interpolation, instead we simply limited interpolation to nearest neighbors defined in terms of the geometry.  This is another issue that we intend to study in detail in the future.  

\bibliography{AMG}{}
\bibliographystyle{plain}

\end{document}